# Multi-objective Optimization Framework for Networked Predictive Controller Design

Sourav Das[a], Saptarshi Das[a,b], and Indranil Pan[a]

a) *Department of Power Engineering, Jadavpur University, Salt Lake Campus, LB-8, Sector 3, Kolkata-700098, India.*

b) *Communications, Signal Processing and Control Group, School of Electronics and Computer Science, University of Southampton, Southampton SO17 1BJ, United Kingdom.*

**Emails:** das_sourav@live.in (S. Das)

saptarshi@pe.jusl.ac.in, s.das@soton.ac.uk (S. Das)

indranil.jj@student.iitd.ac.in, indranil@pe.jusl.ac.in (I. Pan)

**Abstract**

Networked Control Systems (NCSs) often suffer from random packet dropouts which deteriorate overall system's stability and performance. To handle the ill effects of random packet losses in feedback control systems, closed over communication network, a state feedback controller with predictive gains has been designed. To achieve improved performance, an optimization based controller design framework has been proposed in this paper with Linear Matrix Inequality (LMI) constraints, to ensure guaranteed stability. Different conflicting objective functions have been optimized with Non-dominated Sorting Genetic Algorithm-II (NSGA-II). The methodology proposed in this paper not only gives guaranteed closed loop stability in the sense of Lyapunov, even in the presence of random packet losses, but also gives an optimization trade-off between two conflicting time domain control objectives.

## 1. Introduction

Networked control system is a distributed control system wherein the information is exchanged among the system's components such as sensors, actuators, controllers, etc. through a shared real-time communication network. The NCS results in reduced wiring, makes maintenance and analysis of system easier and increases flexibility of architecture, thus making the control system to be cost effective. So now-a-days, NCSs are strongly recommended over conventional point-to-point control systems. NCSs are widely applicable in factory automation, robotics [1]-[2], in safety critical applications and hazardous environments like space excursions, nuclear power plants, motorized vehicles etc. But in NCSs, random packet drop-out is one of the major causes of performance deterioration of the system and often might lead to system instability [3]-[4]. This might happen, inspite of the nominal system being stable under ideal conditions. Thus, stabilization of conventional state-feedback control loops and the same in the presence of random packet-losses should be dealt in a different way.

In Mu *et al.* [5], the classical model predictive control (MPC) scheme has been used to tackle the unknown network induced delay and packet dropouts. In that paper, the receding horizon controller continuously performs an online optimization at every time instant to predict the sequence of control signals based on received data. The predicted control signals are lumped in a packet and transmitted to actuator at the plant side. The delay compensator logic at the actuator end decides the appropriate control signal from the available predicted signal set that can be used to compensate the time delay and packet losses. Here the controller has to finish the computation within a small fraction of sampling instant. This method thus imposes a hard time limit in real time control. In contrast to the above methodology, the present paper does not use the concept of online optimization to tackle packet-losses in NCS. Rather, to reduce the online computation, static state-feedback gains are used which would provide future control signals if the present and subsequent packet is dropped. The concept has been illustrated in Pan *et al.* [6]. In the present predictive NCS scheme, the estimated present and future control signals are lumped in a packet and transmitted through the shared network. There is a buffer which stores the packet before applying the control signal to the plant. The buffer has a pre-programmed logic which decides the control signal to be applied at each sampling instant depending on the previous packet history. When there is no packet loss, only the first estimated control signal is applied to the plant at every sampling instant, and rest of the control signals are discarded. But when the NCS suffers packet loss, some of these packets containing the effective control signals also get dropped. To compensate for these time instants, the buffer applies estimated control signal from last received data packet to the plant, depending upon the receiving instant of the packets. The proposed predictive state feedback controller design should not be confused with the well-known MPC algorithms. In MPC, an iterative method is employed which performs a finite horizon optimization of the plant model at each time step. A numerical optimization algorithm minimizes a user defined control cost at each time step and predicts the next N control signals that must be supplied for the optimal operation of the plant. Only the first control signal from this set of predictive control signals is then used to update the plant state and the whole process is again repeated at subsequent time steps. These algorithms work well for nonlinear, time delayed or other complicated process control plants. In the present case, the general philosophy of predictive control is used to compensate the arbitrary packet loss in the network. The packet loss is a stochastic phenomenon and the designer has no idea about the sequence in which the packet drops would occur, although he knows the number of successive packet drops that can occur in a given network loading scenario. So a generalized formulation needs to be done in which the controller gives guaranteed stability and good performance, irrespective of arbitrary order of packet drops. The term "predictive" refers to a-priori determination of future control actions with an anticipation of stochastic packet loss in the shared communication medium. The control policy can be visualized as a set of linear state feedback controllers in parallel for each case of the different packet drop scenarios. The control action from each of the controller is calculated and encoded in a packet and sent over the network. Depending on the arrival of the packets and the drop outs, the logic at the buffer decides which control action to give to the actuator. This buffer logic mechanism along with the set of parallel linear controllers is what performs the predictive strategy. The advantage of this method over the other MPC methods is that there is no need of a finite horizon optimization (and hence increased computational cost) at each time step since the control design method is offline and the prediction occurs by the synergistic logic and interaction between the multiple set of linear controller and the buffer logic. Also this can be implemented by simple state feedback controllers as opposed to more complex MPC paradigms.

The Asynchronous Dynamical System (ADS) is another approach of modeling the NCS with packet loss, which is used in Zhang *et al.* [4]. But ADS formulations involve Bilinear Matrix Inequality (BMI) constraints and there is no standard solution technique for BMIs constraints in polynomial time. In most cases the BMI problems are solved using random search techniques. For this reason, the ADS approach has not been considered in the present work. In this paper, the NCS is modeled as a set of switched Linear Time Invariant (LTI) systems with arbitrary switching between them caused by random data loss. This mechanism avoids the problems in ADSs. In switched system scheme, the equations are formulated using linear matrix inequality (LMI) constraints which can be solved using standard off-the-shelf LMI solvers to ensure Lyapunov stability for the system, even in the presence of packet-losses. In the present approach, the LMI based stability criteria is solved as a sub-problem in each iteration and the controller design is done by global optimization technique using Multi-Objective Genetic Algorithm (MOGA) [7].

In Xiong and Lam [8], the switched system methodology is used for stabilization of NCSs but the performance criteria are not considered in the design problem. But for industrial uses such as networked process control applications, the quality of control is also a major consideration and not mere stabilization. With the present formulation, guaranteed Lyapunov stability and optimum time domain performance both can be achieved in NCSs in the presence of random packet losses. However, in a practical NCS, often conflicting objectives are encountered which needs to be optimized to get a trade-off between them. Contradictory objectives imply that if the performance of one objective is improved, performance of other objective gets degraded and it is not possible to minimize both the objectives simultaneously.

In Pan *et al.* [6], Particle Swarm Optimization (PSO) has been used to maximize the stabilizing region of the predictive networked controller. From their simulations it is clear that the stabilizing region is patchy and discontinuous and consequently cannot be visualized as a standard continuous nonlinear objective function. Some of the continuous nonlinear objective functions can also be solved using standard LMI with special transformation techniques like the Schur complement, change of variables etc [9]. On the other hand evolutionary and swarm optimization algorithms have been found to be expedient in optimizing functions which are discontinuous, noisy or dynamically varying which often occur in standard control engineering problems [10]. Hence in the present approach, the controller design problem is attempted through population based multi-objective evolutionary optimization techniques. LMI based techniques are preferred over evolutionary approaches, in cases where the problem is convex, since they can be solved in polynomial time and have the feature of guaranteed convergence. However in the present case, the objective functions are non-convex due to their discontinuous nature. To further explain this point, a reference can be made to the stability region in the state feedback controller parameter space as in Pan *et al.* [6] which is clearly a discontinuous function. As a subset of this problem, the present paper seeks the optimum performance inside the stable regions only. Hence the region, in which the multiobjective optimization searches needs to be confined, is also discontinuous in nature in order to achieve optimum time domain performance. This is due to the fact that only stabilizing controller gains can produce finite integral performance index associated with them where the multi-objective search is to be carried out. Thus, even though the evolutionary techniques have the issues of consistency and guaranteed convergence, these have been resorted to and LMI techniques have been used as a sub-problem in the present scenario.

Multiobjective optimization also can be done using LMI [11]-[12]. However, the functions need to be convex for this type of optimization which is not valid for the present networked predictive controller design scenario. So MOGA has been used for simultaneous minimization of multiple contradictory objectives for this non-convex, discontinuous objective function. In this paper, the multi-objective NSGA-II is used to find out the Pareto optimal solutions for the controller parameters as a trade-off between different contradictory performance objectives [13]-[20]. MOGA provides a powerful population-based search and a flexible optimal control system can be designed with the user's choice of putting weights of those objectives [16]-[20]. It is used to optimize the objective functions, represented by few time domain performance indices for those gains which satisfy the LMIs for state-feedback control with packet-losses.

The rest of the paper is organised as follows. In section 2 the NCS is analysed as a switched system and the LMI formulation for the predictive controller is outlined. The NSGA II algorithm is also introduced and various contradictory objective functions are discussed. An illustrative example is given in section 3, along with credible simulation results to demonstrate the proposed methodology. The paper ends in section 4 with the conclusions followed by the references.

*Notations:*

Throughout this paper, $\mathbb{R}^n$ and $\mathbb{R}^{n \times m}$ denote the *n* dimensional Euclidean space and the set of all $n \times m$ real matrices respectively. $\|\cdot\|$ refers to the Euclidean norm for vectors and induced 2-norm for matrices. $I$ and $\phi$ are the identity and null matrix respectively with appropriate dimensions. $\mathbb{Z}^+$ is the set of positive integers. The superscript "$T$" indicates matrix transpose operation.

## 2. Theoretical formulation

### 2.1 LMI formulation for the state feedback NCS with predictive gains considering drop in both the forward and the feedback channels

The model of the plant concerned in this paper is shown in the Fig. 1. The discrete time state space model of the process is represented as:

$$x(k+1) = Fx(k) + Gu(k) \tag{1}$$

where, $k \in \mathbb{Z}_+$ is the time step, $x(k) \in \mathbb{R}^n$ is the system state and $u(k) \in \mathbb{R}^m$ is the control input; $x_0 := x(0) \in \mathbb{R}^n$ is the initial state. $F$ and $G$ are two system matrices of appropriate dimensions.

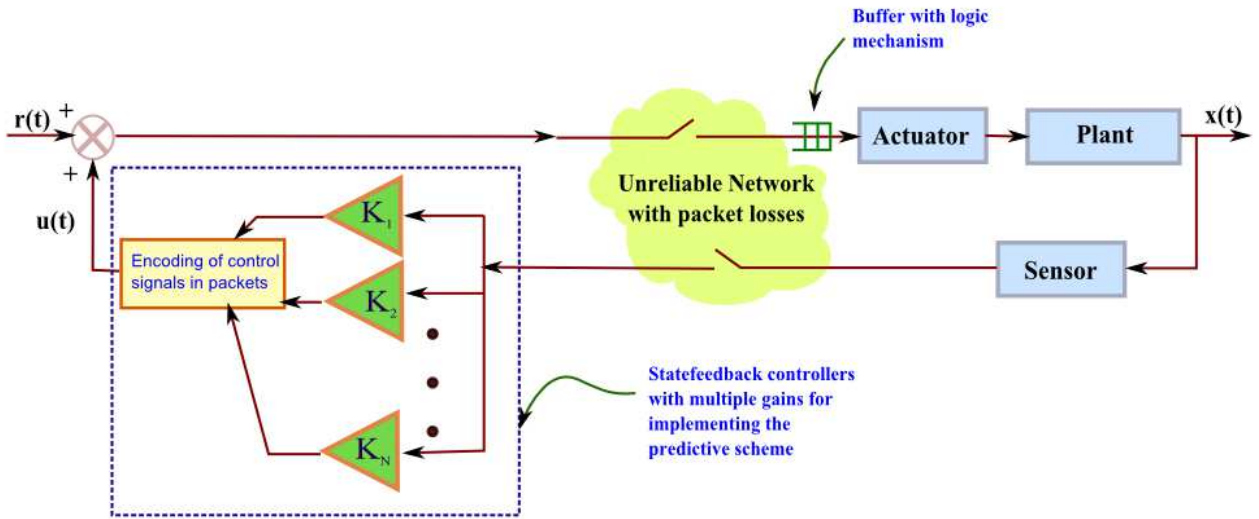

Fig. 1. Schematic diagram of considered predictive NCS.

Fig. 1 shows the schematic of the NCS with the predictive feedback controller and the unreliable network having packet drops in the forward and the feedback paths. The sensor is assumed to be time driven. It samples the plant at fixed multiples of the sampling time and sends the data to the controller in the form of packets along with a time stamp. The controller is a full state feedback controller and full state information of the plant is assumed to be available. The controller can be time-driven or event-driven. But in the peak load situation of the network, a number of packets can arrive within the same sampling instant for the event driven case, which is undesirable for safety critical applications. For this reason, in the present paper the controller is chosen as time-driven. The controller encapsulates the present and next predicted control signals into a data packet and sends it over the network. The control signal is given by

$$u = K_z x \qquad (2)$$

where, $K_z$ is the feedback gain to be designed and $z \in \{1, 2, ..., M_{drop}\}$.

The buffer situated on the plant side has a logic element in it. It receives the control packet and applies the proper control signal to the actuator based on the arrival history of packets. The logic of the buffer is explained in Definition 1. The time lines of the system sensor, the controller and the actuator are assumed to be synchronized with a skew to take care of the inherent delays. Thus if $\tau_1$ be the time skew between the sensor and the controller, $\tau_2$ be the time skew between the controller and the actuator, then the following is assumed to hold: $T_s > \tau_1 + \tau_2$, where $T_s$ is the sampling time of the system. Thus a packet with a particular time step is considered to be dropped if it does not reach the actuator in the same time step. Any out of order packet is discarded by the buffer and considered to be dropped. Let $S := \{i_1, i_2, ...\} \subseteq \mathbb{Z}^+$ indicate the sequence of time points of the effective packets which are transmitted successfully from the sensor to the actuator buffer (Fig. 2). Based on above notations, packet loss process is defined as

$$\xi(i_m) := \{i_{m+1} - i_m \mid i_m \in S\} \qquad (3)$$

which means from $i_m$ to $i_{m+1}$, maximum $\xi(i_m) - 1$ number of packets can be dropped consecutively. Also let, $M_{drop} := \max_{i_m \in S}\{\xi(i_m)\}$. Hence, $\xi(i_m)$ can take values in the finite set $\mu := \{1, 2, ..., M_{drop}\}$. The packet loss process is assumed to take random values in $\mu$. At every sampling instant, the controller sends packets encapsulated with present and future control signal. The control signal can be represented as $u_{pq} = K_q x(i_p)$ where $K_q = [K_{q1}\ K_{q2}\ \cdots\ K_{qn}]$, $p \in \{1, 2, 3, ...\}$ denotes the sample time of effective packets and $q \in \{1, 2, 3 ..., M_{drop}\}$ represents the index numbers of the next $M_{drop}$ samples of control signals from the present sampling instant.

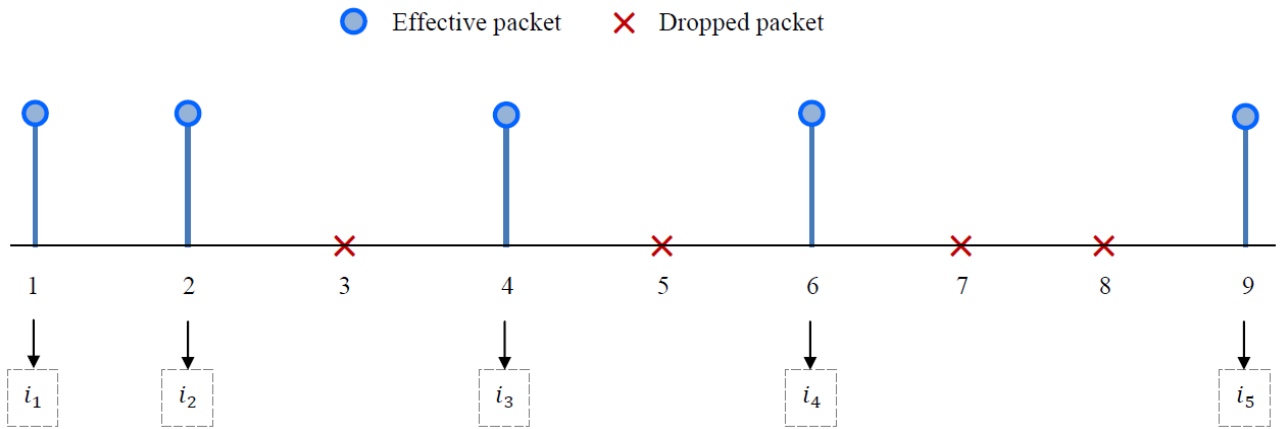

Fig. 2. Illustration of packet drops in the predictive NCS scheme.

As for example, consider a case where, maximum two consecutive packets can be dropped i.e. $M_{drop} = 3$, The situation of packet drops in NCS can be as shown in Fig.2. It means that 1st, 2nd, 4th, 6th, 9th sensor packets are effective packets used to control the plant. Whereas the 3rd, 5th, 7th and $8^{th}$ sensor packets are considered as dropped and they are not used. At every $k^{th}$ $(k=1,2,3,\cdots,9)$ sampling instant the control input used to control the plant is given by $\{u_{11},u_{21},u_{22},u_{41},u_{42},u_{61},u_{62},u_{63},u_{91}\}$ respectively.

***Definition 1 [6]:*** The buffer logic which decides the control signal $u$ to the actuator at each sampling instant $k$, is defined as $u(k):=\psi_{k-\rho+1}(\rho)$ where $\rho=k-i_m+1$. Here, $\psi$ denotes the buffer, the subscript $k-\rho+1$ denotes the time index the buffer was last updated by an effective packet. $\rho$ denotes the array index of the buffer from which the control signal is applied to the actuator. The plant (1) at each instant is updated by the control values in the buffer depending on $\xi(i_m)$. Thus at any given instant the plant can be controlled by a packet index $\psi\in\mu$ in the buffer.

Buffer states corresponding to the effective packets as in Fig. 2, are shown in Table 1, The buffer is a memory device which will register maximum three data at every sampling instant in an array for $M_{drop}=3$. The data of first column is used to update the plant at every instant which is shown in the column "Present" in Table 1. Status of next predicted data is shown in columns "1$^{st}$ Predicted Sample" and "2$^{nd}$ Predicted Sample" respectively. At sampling instants 1 and 2, the buffer is updated with the data of two consecutive successfully transmitted packets containing control signals $\{u_{11},u_{12},u_{13}\}$ and $\{u_{21},u_{22},u_{23}\}$ respectively. When the buffer sends the current control signals such as $u_{11}$ and $u_{21}$ for the instants 1 and 2 respectively, other control signals followed by those signals, predicted to compensate for next instants, are shifted left by one. At instant 3, the packet containing current control signal from controller fails to reach the buffer. So the predicted data $u_{22}$ which is already transmitted to the buffer with the previous packet is used to update the plant and $u_{23}$ is shifted left. At next instant, the packet is again lost. Hence $u_{23}$ is sent to the plant. At instant 4, the packet successfully reaches to the buffer and update it with the control signal data $\{u_{41},u_{42},u_{43}\}$. In this way the buffer is being updated for next instants also.

Table 1

Buffer states corresponding to packet arrival in Fig. 2

| Sampling Instant | Present | 1$^{st}$ Predicted Sample | 2$^{nd}$ Predicted Sample |
|---|---|---|---|
| 1 | $u_{11}$ | $u_{12}$ | $u_{13}$ |
| 2 | $u_{21}$ | $u_{22}$ | $u_{23}$ |
| 3 | $u_{22}$ | $u_{23}$ | - |
| 4 | $u_{41}$ | $u_{42}$ | $u_{43}$ |
| 5 | $u_{42}$ | $u_{43}$ | - |
| 6 | $u_{61}$ | $u_{62}$ | $u_{63}$ |
| 7 | $u_{62}$ | $u_{63}$ | - |
| 8 | $u_{63}$ | - | - |
| 9 | $u_{91}$ | $u_{92}$ | $u_{93}$ |

*2.2 Discrete time switched system stability of NCS using predictive controller*

This switched system approach essentially incorporates the idea of the lifted sampling period. Consider that every alternate packet of the control loop is dropped. This is equivalent to sampling it at twice the sampling time. For two consecutive packet drops the equivalence is with thrice the sampling time. In other words, the sampling timeline is "lifted" to a new timeline which is some integer multiple of the base sampling rate. In the present analysis the NCS is modelled in such a manner that the packet drop process is modelled as switching between these different lifted timelines. In the NCS modelling scenario the switching between the finite number of states is arbitrary and hence there is no control over the order of the switching. Hence this is not a stabilization problem where a switching schedule can be designed to stabilize the system. Asymptotic stability of the overall system must be guaranteed in the presence of arbitrary switching between the sub-systems.

Let us consider the augmented state vector $\Gamma(k) = \left[ x(k)\ x(k-1) \ldots x(k-M_{drop}+1) \right]^T$. The NCS with the state feedback predictive controllers can be cast in the form of a discrete time switched system given by:

$$\Gamma(k+1) = \Phi_{\sigma(k)} \Gamma(k) \tag{4}$$

where, $\sigma(k)$ is a piecewise constant function, known as a switched signal, which takes values in the finite set $\Lambda := \{1, 2, \ldots, M_{drop}\}$. $\Phi_{\sigma(k)} \in \mathbb{R}^{nM_{drop} \times nM_{drop}}$ is of the generalized form:

$$\Phi_{\sigma(k)} = \begin{bmatrix} \overline{F}_{n \times n} & H_{n \times n(M_{drop}-1)} \\ I_{n(M_{drop}-1) \times n(M_{drop}-1)} & \phi_{n(M_{drop}-1) \times n} \end{bmatrix} \tag{5}$$

where, $\overline{F} = \begin{cases} F + GK_1 & \text{for } \sigma(k) = 1 \\ F & \forall \sigma(k) \in \{2, \cdots, M_{drop}\} \end{cases}$;

$H(1, \sigma(k)) = GK_{\sigma(k)}\ \forall \sigma(k) \in \{2, \cdots, M_{drop}\}$, and $H(1, t) = \phi_{n \times n}\ \forall t \neq \sigma(k)$

***Theorem 1 [6]:*** The NCS defined by the switched system (4) is asymptotically stable for the arbitrary packet loss process $\xi(i_m)$, if for $i, j \in \Lambda$, there exists positive definite matrices, $P_i$ and $P_j$ satisfying the following set of LMIs.

$$\Phi_i^T P_j \Phi_i - P_i < 0 \tag{6}$$

where, $\Phi_i$ is of the form (5) with specified controller gains $K_l\ \forall l \in \Lambda$.

***Proof:*** For the switched system (4) let us define multiple quadratic Lyapunov functions of the form (7), for each switched state.

$$V(k) = \Gamma^T(k) P_{\sigma(k)} \Gamma(k) \tag{7}$$

where, $P_{\sigma(k)}$ are symmetric positive definite matrices $\forall\ \sigma(k) \in \Lambda$.

Let the value of $\sigma(k)$ at the $k^{th}$ and $(k+1)^{th}$ time instant be $i$ and $j$ respectively, where $i, j \in \Lambda$. The difference of the Lyapunov function between the two instants of time is given by:

$$\begin{aligned} \Delta V(k) &= V(k+1) - V(k) \\ &= \Gamma^T(k+1) P_j \Gamma(k+1) - \Gamma^T(k) P_i \Gamma(k) \\ &= \Gamma^T(k) \left( \Phi_i^T P_j \Phi_i - P_i \right) \Gamma(k) \end{aligned} \tag{8}$$

For any $\Gamma(k) \neq 0$, $\Delta V(k) < 0$ if (6) holds.

Thus, $\lim_{k \to \infty} V(k) = 0 \ \forall \ k \in \mathbb{Z}_+$. Hence the system (4) is asymptotically stable. □

Therefore, the linear equations (6) are solvable if matrices $\Phi_i$ are stable i.e. all eigen-values of $\Phi_i$ should have negative real part. In other word, we have to show that the LMIs are feasible for some value of $P_i > 0$ and $P_j > 0$.

*2.3 MOGA based controller design with optimized performance and guaranteed stability*

Multi Objective Genetic Algorithm (MOGA) is a population based optimization algorithm based on Darwinian principle of survival of the fittest. This algorithm transforms a set of solution variables into population of solutions depending upon individual fitness value through reproduction. Reproduction implies that solution vectors with higher fitness values can produce more replicas of themselves in the next generation. Usually a parameter called the elite count is used which represents the number of fittest solution vectors that will definitely go to the next generation. But increasing the elite count may result in domination of the fitter individuals obtained earlier in the simulation process and as such will result in less effective solutions. Hence this parameter is generally a small fraction of the total population size. After the selection of the fittest solutions, different operations like crossover and mutation take place among them to produce better and effective solutions in each generation. In crossover, information is exchanged between a pair of solution vectors so that good solutions can frequently arrive in the next generation. It makes the search process to converge towards overall best solutions with similar characteristics. In mutation any randomly selected portion of solution vector is altered, introducing diversity to similar solutions. It helps to avoid local minima and initiates to search the unexplored regions of Pareto-front to find new non-dominated set of solutions. It randomly changes the information of individual solution.

The multiobjective optimization problems can be solved by transforming the objectives into single objective by assigning weights to individual objectives and any of single objective optimization algorithms can be applied to solve the problem. However in this case the solution will depend upon the specific choice of the weights. Sometimes appropriate choice of these weights can be very difficult to guess for the designer. Even a small change in the weights can produce very different solutions. On the other hand MOGA, being a population based search algorithm, can simultaneously search different regions of the solution space. This property is helpful to find the different set of solutions for discontinuous or non-convex problem of switched systems. MOGA can search a set of solutions keeping one objective under acceptable level without being dominated by other objectives. Thus simultaneous optimization of multiple objectives can be achieved. For these numerous advantages, MOGA is chosen for the optimization of conflicting objectives in the present design procedure.

In this paper, NSGA-II algorithm is used for multi-objective optimization [21]. The crossover and mutation operations are similar to those in single objective genetic algorithm excluding the selection procedure. Before the selection is performed, the non-dominated solution of the population is given a high dummy fitness value. Then these solutions are ignored temporarily and a new non-dominated set of solutions are formed from rest of the population and is assigned with lesser fitness value. This process continues until all the solutions are assigned with a fitness value. After that all solutions are reproduced according to individual fitness value. As the individuals in the first front have highest fitness value, they produce more replicas than others. This makes them to converge faster on the non-dominated space. A generalized multi-objective optimization framework can be defined as follows:

Minimize $F'(x) = (f_1(x), f_2(x), ..., f_m(x))$ such that $x \in \Omega$ (9)

where $\Omega$ is the decision space, $\mathbb{R}^m$ is the objective space, and $F': \Omega \rightarrow \mathbb{R}^m$ consists of $m$ real valued objective functions.

Let $v = \{v_1, ..., v_m\}$, $w = \{w_1, ..., w_m\} \in \mathbb{R}^m$ be two vectors. Here, $v$ is said to dominate $w$ if $v_i < w_i \ \forall \ i \in \{1, 2, ..., m\}$ and $v \neq w$. A point $x^* \in \Omega$ is called Pareto optimal if $\nexists \ x \mid x \in \Omega$ such that $F'(x)$ dominates $F'(x^*)$. The set of all Pareto optimal points, denoted by PS is called the Pareto set. The set of all Pareto objective vectors, $PF = \{F'(x) \in \mathbb{R}^m, x \in PS\}$, is called the Pareto Front. This implies that no other feasible objective vector exists which can improve one objective function without simultaneous worsening of some other objective function.

Multi-objective Evolutionary Algorithms (MOEAs) which use non-dominated sorting and sharing have higher computational complexity, uses a non-elitist approach and requires the specification of a sharing parameter. The non-dominated sorting genetic algorithm (NSGA-II), removes these problems and is able to find a better spread of solutions and better convergence near the actual Pareto optimal front [21]. The pseudo code for the NSGA-II is as shown below [21]-[22].

---

NSGA II Algorithm

Step 1: generate population $Y_0$ randomly

Step 2: set $Y_0 = (F_1', F_2', ...) = $ non-dominated-sort$(Y_0)$

Step 3: for all $F_i' \in Y_0$

    crowding-distance-assignment$(F_i')$

Step 4: set t=0

    while (not completed)

        generate child population $Q_t$ from $Y_t$

        set $R_t = Y_t \cup Q_t$

        set $F' = (F_1', F_2', ...) = $ non-dominated-sort$(R_t)$

        set $Y_{t+1} = \phi$

        i=1

        while $|Y_{t+1}| + |F_i'| < N$

            crowding-distance-assignment$(F_i')$

            $Y_{t+1} = Y_{t+1} \cup F_i'$

            i=i+1

        end

        sort $F_i'$ on crowding distances

        set $Y_{t+1} = Y_{t+1} \cup F_i'[1:(N - |Y_{t+1}|)]$

        set $t = t + 1$

    end

return $F_1'$

---

The NSGA II algorithm converts $M$ different objectives into one fitness measure by composing distinct fronts which are sorted based on the principle of non-domination. In the process of fitness assignment, the solution set not dominated by any other solutions in the population is designated as the first front $F_1'$ and the solutions are given the highest fitness value. These solutions are then excluded and the second non dominated front from the remaining population $F_2'$ is created and ascribed the second highest fitness. This method is iterated until all the solutions are assigned a fitness value. Crowding distances are the normalized distances between a solution vector and its closest neighbouring solution vectors in each of the fronts. All the constituent elements of the front are assigned crowding distances to be later used for niching. The selection is achieved in tournaments of size 2 according to the following logic.

a) If the solution vector lies on a lower front than its opponent, then it is selected.
b) If both the solution vectors are on the same front, then the solution with the highest crowding distance wins. This is done to retain the solution vectors in those regions of the front which are scarcely populated.

## 3. Formulation and simulation of the conflicting objective functions

Different methods of controller design are described now by optimizing two different set of objective functions simultaneously. As shown in Pan *et al.* [6], time domain performance decreases with increase in stability region for the predictive control. In control engineering, along with good relative stability there are also several other performance requirements which should be met in order to design a controller. These can be small steady state error, good transient response like small peak time, rise time, small overshoot, and robustness of parameters. All the performance parameters are inter-related to each other and if one parameter is optimized, other parameters will also change automatically. So there is a requirement of optimizing more than one parameter simultaneously, rather than one parameter at a time. Hence to achieve this, three design techniques are introduced, each consisting of two contradictory objective functions which are traded-off to achieve better performance.

### 3.1. First design trade-off for controller tuning

The first controller design comprises of two conflicting cost functions given by (10) and (11)

$$J_1 = \sum_{l=1}^{n} \sum_{k=0}^{\infty} (kT_s) |x_l(k)| \tag{10}$$

and $J_2 = \sum_{l=1}^{n} \sum_{k=0}^{\infty} u_l^2(k)$ (11)

where $n$ denotes the number of state variables of the system; $x_l(k)$ is the $l^{th}$ state at $k^{th}$ time instant for the system given by equation (1).

The performance index $J_1$ ensures fast transient response or fast settling of the system's states. Here, time multiplication term with the state variable minimizes the oscillation of the state as time increases and absolute value of state variable minimizes the percentage overshoot of the respective states. But on the other hand control signal is equally important for achieving this specific task. The required control signal should not be allowed to increase indefinitely as it can lead to saturation of the actuator. In order to prevent saturation of the actuator and also to reduce the size of the actuator, squared control signal $J_2$ is taken as another cost function. Squared value of the control signal is taken to put extra penalty on the higher values. Hence by balancing the cost functions $J_1$ and $J_2$ we can reduce oscillations of the state variable, reduce its peak overshoot, settling time while using the minimum possible control signal, due to the multi-objective optimization based design framework of the problem. During the simulations, the clock output is discretized and multiplied with each state at every sampling instant. These results are added to obtain the numerical values of objective $J_1$ and $J_2$, in each iteration of the MOGA.

### 3.2. Second design trade-off for controller tuning

Another objective function $J_3$ is given by (12) which enforces a smooth variation of the state trajectories, even in the presence of predictive gains to compensate for the dropped data packets.

$$J_3 = \sum_{l=1}^{n} \sum_{k=0}^{\infty} \|x_l(k) - \hat{x}_l(k)\|$$
$$\hat{x}_l = smooth(x_l)$$
(12)

Here, function "smooth" reduces the small local oscillations in the time series of state variables ($x_l$) using a moving average (MA) filter within each function evaluation of MOGA. Due to the introduction of predictive gains so as to maximize stabilizing region in Pan *et al.* [6] for packet drops in NCSs, such oscillations occurs in the state-variables and smoothing of the state excursions becomes a necessity. The smoothing results are stored in $\hat{x}_l$ and its difference with the original state vector can be considered as the error and its Euclidean distance as the error index which needs to be minimized. The moving average filter smoothes the time series of the state variables by replacing each data point with the average of the neighbouring data points defined within a given span. The span used for the moving average filter consists of consecutive 5 samples. At every sampling instant, the value of each state is stored in a buffer. Then those values are passed through the 'smooth' function of MATLAB to obtain local oscillation free data set $\hat{x}_l$. The 2-norm of difference between original data and smoothed data of all states has been considered as the third objective function $J_3$. Pan *et al.* [6] have shown that state feedback controller with predictive gains with high robust stability yields jittery nature of the state excursions. So the objective function $J_3$ is chosen to minimize jittery effect of the states which is not desirable for many practical control applications. But in order to reduce the jittery effect, much more control signal will be needed, which may again saturate the actuator. So, integral squared control signal $J_2$ in discrete form has again been taken as the other objective function to study the second design trade-off.

### 3.3. Third design trade-off for controller tuning

It is known that if there is an attempt to make the settling time short, as its consequence the peak value of the state vector increases and vice-versa. For small settling time, the state feedback controller gains should be increased and consequently the controller output also increases. But this will lead to increase in overshoot of the state variables. Hence here these two conflicting features have been used in the multi-objective optimization framework to design the predictive state feedback controllers as also studied in [23] for single objective optimization. Here the performance indices $J_4$ (13) and $J_5$ (14) impose penalties on high peak value and long settling time respectively.

$$J_4 = \sum_{l=1}^{n} \left| \frac{M_{peak}(x_l)}{x_l^0} \right| \qquad (13)$$

where $x_l^0$ is the initial value of $l^{th}$ state vector. Peak value of $l^{th}$ state vector is denoted by $M_{peak}(x_l)$. Since, in many practical applications initial value of the state variables ($x_l^0$) may be different, therefore they are normalized in (13) with respect to $x_l^0$ to ensure that overshoots are considered in the same normalized scale in the optimization process. It is to be noted that in the evaluation of objective function (13), zero initial states have to be replaced with small initial values for to avoid division by zero.

$$J_5 = \sum_{l=1}^{n} ST_l \qquad (14)$$

where $ST_l$ is the settling time of $l^{th}$ state vector.

In the simulation of performance index (14), the absolute values of a state are stored in a buffer at every sampling time. If maximum value of last 10 consecutive samples becomes less than a certain tolerance (taken as ±0.02 or 2% criterion in this case), it is confirmed that the respective state response reaches inside the tolerance band of its final value and that instant can be considered as the settling time (ST).

Here controller gains are randomly generated using MOGA so as to study design trade-off between above control objectives. Within the multi-objective optimization process, it is first checked using the LMI Toolbox "YALMIP" [24] whether the randomly generated controller gains satisfy the LMIs as given in equation (6) derived from Lyapunov stability criteria. In case the stability criteria is not satisfied, a high value of objective function is assigned which automatically makes the solution inferior and the algorithm steers off from these regions of the search space in the consequent generations. The objective functions mentioned above are stochastic in nature because it is not possible to predict at which time instant the packet losses occur. Hence both the objective functions in multi-objective optimization for every set of controller design have been simulated multiple times using same set of solution vectors (controller gains) and average values of objective functions are obtained. In the simulation, a random variable is used to decide whether the packet will be dropped or not. The random variable generates random value between 0 and 1. If the generated value is less than 0.8, the packet is considered as drop. Otherwise the packet is taken as effective packet.

Similar LMI formulation for networked process control application can be found in recent literature e.g. using dynamic matrix control [25], observer based $H_\infty$ control [26] for nonlinear systems, receding horizon $H_\infty$ control [27], robust fault detection filter [28], distributed model predictive control [29], guaranteed cost MIMO control [30], event-triggered control [31] etc. The proposed methodology searches for Pareto optimal solutions for controller gains with different time domain control objectives with the LMI criterion serving as an inherent checking condition for guaranteed stability. The flowchart of a systematic closed loop networked controller design has been shown in Fig. 3.

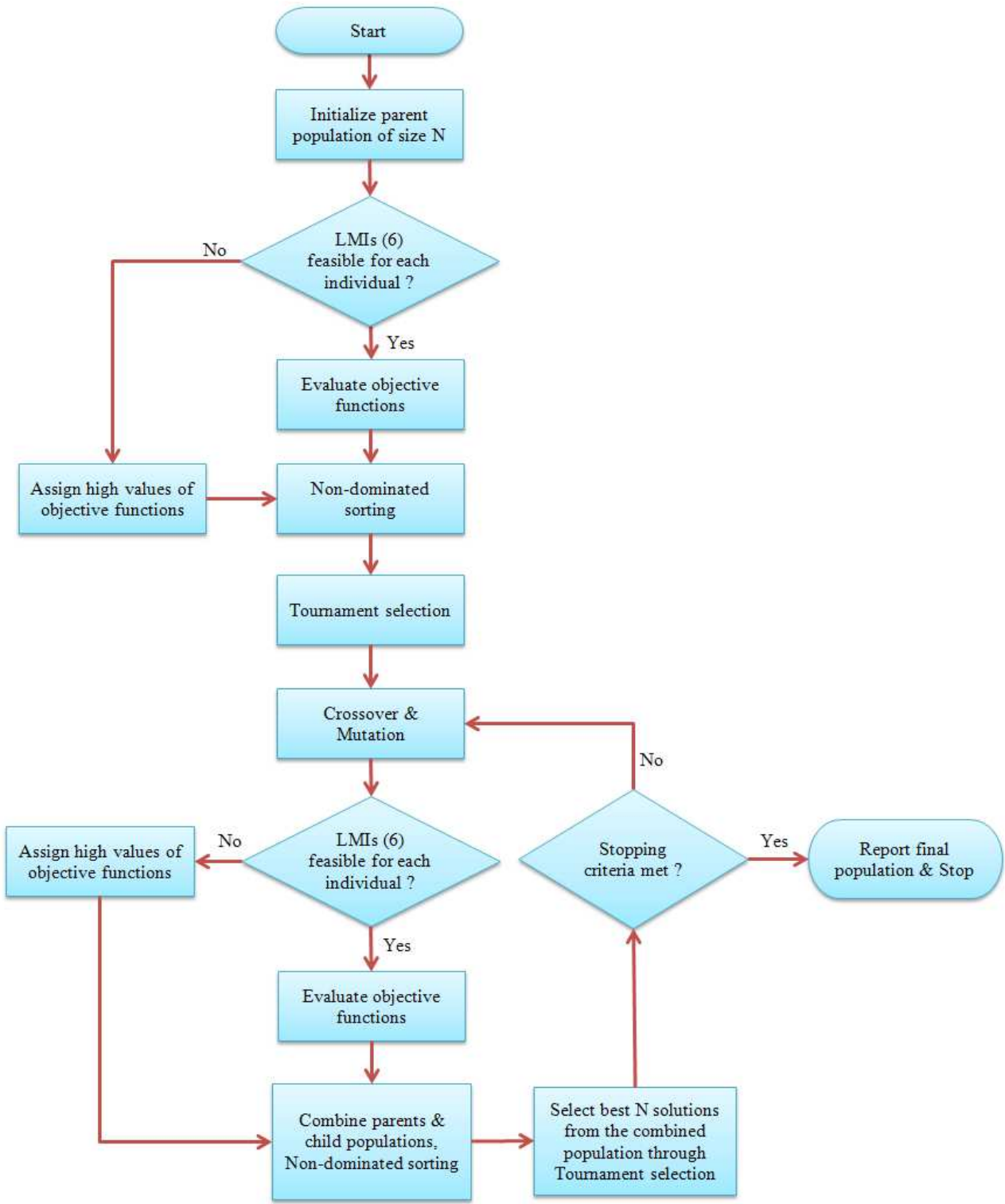

Fig. 3. Flowchart for the proposed controller design technique.

## 4. Illustrative simulation examples

### *4.1. DC motor plant*

Let us consider, an NCS with a continuous time dc motor plant (15) as in [32]-[34] with the state variables considered as the angular position and angular velocity. The dc motor parameters are detailed in [34]. The continuous time plant is unstable with a pole very close to the origin and the other pole being non-dominating is far away from the imaginary axis of the complex s-plane.

$$\dot{x}(t) = \begin{bmatrix} 0 & 1 \\ 1 & -217.4 \end{bmatrix} x(t) + \begin{bmatrix} 0 \\ 1669.5 \end{bmatrix} u(t) \quad (15)$$

By taking sampling period as 0.05s, we can get discrete time plant of the structure (1) having the system matrices given by (16)

$$F = \begin{bmatrix} 1.0002 & 0.0046 \\ 0.0046 & 0 \end{bmatrix}, \quad G = \begin{bmatrix} 0.3487 \\ 7.6807 \end{bmatrix} \quad (16)$$

Initial values of the state variables are chosen as $x_0 = [3 \ -2]^T$. $M_{drop} = 3$ is assumed to illustrate the proposed design technique. As $M_{drop} = 3$, the augmented state matrix $\Gamma(k) = [x(k) \ x(k-1) \ x(k-2)]^T$ is introduced into NCS. Then according to equation (4) the closed loop NCS can be expressed as

$$\Gamma(k+1) = \Phi_{\sigma(k)} \Gamma(k) \text{ for } \sigma(k) \in \{1,2,3\} \quad (17)$$

where, $\Phi_1 = \begin{bmatrix} F+GK_1 & \phi & \phi \\ I & \phi & \phi \\ \phi & I & \phi \end{bmatrix}$, $\Phi_2 = \begin{bmatrix} F & GK_2 & \phi \\ I & \phi & \phi \\ \phi & I & \phi \end{bmatrix}$ and $\Phi_3 = \begin{bmatrix} F & \phi & GK_3 \\ I & \phi & \phi \\ \phi & I & \phi \end{bmatrix}$.

For the NSGA-II simulations the population size is set as 90, cross-over fraction is chosen as to be 0.8 and Pareto-front population factor as 0.35. Mutation fraction is considered as to be 0.2. The algorithm is run up to a maximum 200 number of generations. The Pareto-front obtained after multi-objective optimization of $J_1$ and $J_2$ is shown in Fig. 4.

From the Pareto-front of Fig. 4, two extreme solutions (solution $A_1$ and $C_1$) and one solution in between (solution $B_1$, which is the median solution) have been chosen for demonstrating their time domain characteristics. The corresponding controller gains are given in Table 2. As the state feedback control loop behaves like a switched system in presence of packet losses, positive state feedback can also stabilize the plant. Hence gains, as assumed in equation (2), can also be positive.

Table 2:

Multi-objective optimization results amongst the objective functions for the dc motor plant

| Design trade-off amongst objective functions | Solution points | Objective functions | | | | | State feedback controller gains | | | | | |
|---|---|---|---|---|---|---|---|---|---|---|---|---|
| | | $J_1$ | $J_2$ | $J_3$ | $J_4$ | $J_5$ | $K_{11}$ | $K_{12}$ | $K_{21}$ | $K_{22}$ | $K_{31}$ | $K_{32}$ |
| $J_1$ and $J_2$ | $A_1$ | 5.393 | 0.090 | - | - | - | -0.155 | 0.003 | -0.047 | 0.036 | -0.152 | -0.040 |
| | $B_1$ | 33.123 | 0.028 | - | - | - | -0.040 | 0.021 | -0.035 | 0.010 | -0.065 | -0.042 |
| | $C_1$ | 62.284 | 0.019 | - | - | - | -0.026 | 0.013 | -0.038 | -0.001 | -0.044 | -0.045 |
| $J_3$ and $J_2$ | $A_2$ | - | 0.056 | 2.000 | - | - | -0.104 | -0.015 | -0.100 | -0.012 | -0.116 | -0.039 |
| | $B_2$ | - | 0.040 | 2.076 | - | - | -0.071 | -0.009 | -0.077 | -0.021 | -0.087 | -0.049 |
| | $C_2$ | - | 0.032 | 2.524 | - | - | -0.058 | -0.007 | -0.058 | -0.035 | -0.080 | -0.048 |
| $J_4$ and $J_5$ | $A_3$ | - | - | - | 2.000 | 13.950 | -0.091 | -0.009 | -0.074 | 0.022 | -0.127 | -0.051 |
| | $B_3$ | - | - | - | 2.406 | 9.350 | -0.116 | 0.026 | -0.054 | 0.075 | -0.135 | -0.012 |
| | $C_3$ | - | - | - | 3.268 | 7.050 | -0.178 | 0.029 | -0.069 | 0.075 | -0.120 | -0.012 |

After getting the Pareto solutions, the simulations have been run for those set of optimal controller gains corresponding to $A_1$, $B_1$ and $C_1$. As the system is stochastic, the values of objectives will not be exactly same as the values which are obtained from the Pareto-front. So the simulation is run for several times and only those time response plots are shown for which values of objectives become most near to the obtained Pareto solution. The time evolution of the state variables and the associated control signals have been depicted in Fig. 5-6 for the three chosen solutions given in Table 2. It is clear from Fig. 5 that the states for the solution $A_1$ settles faster than others, whereas states for $C_1$ take more time to settle down to the final value. So in case of $A_1$, more control action is needed to settle down the states in much shorter time. On the other hand, to settle the states within greater time, control action needed is less for solution $C_1$ as shown in Fig. 6. These situations justify the design trade-off between the two chosen performance criteria.

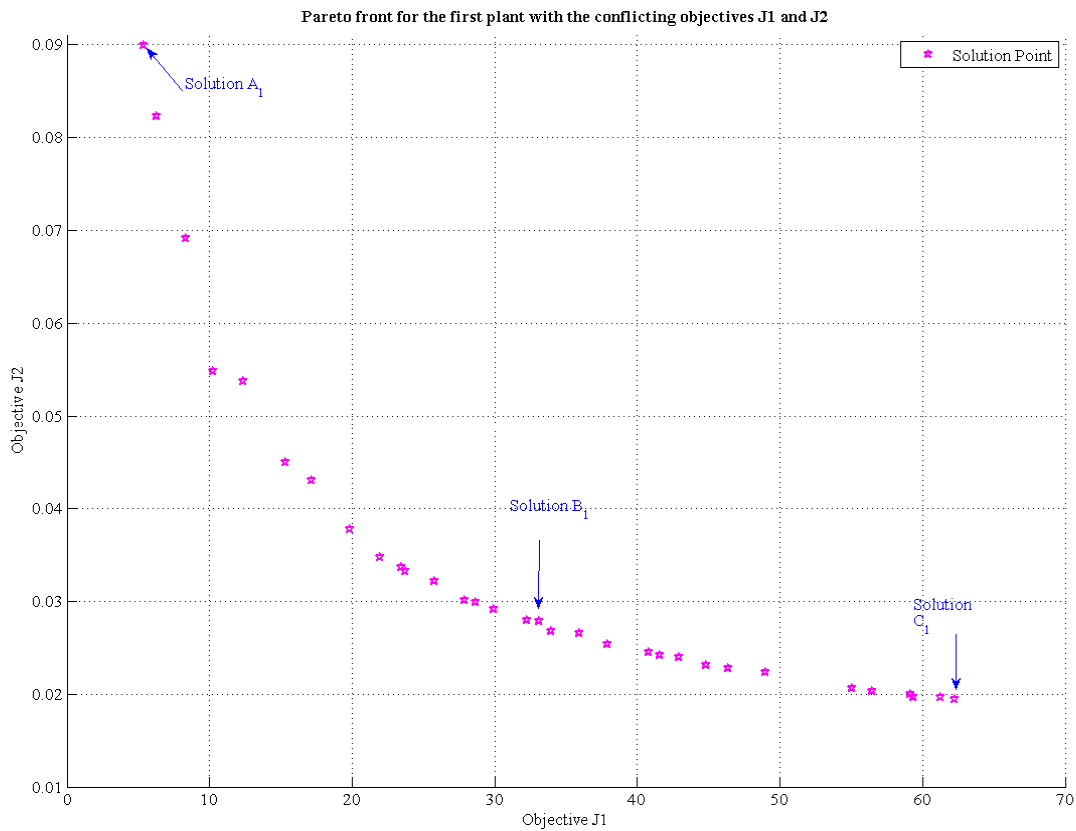

Fig. 4. Pareto optimal front for the conflicting objective functions $J_1$ and $J_2$

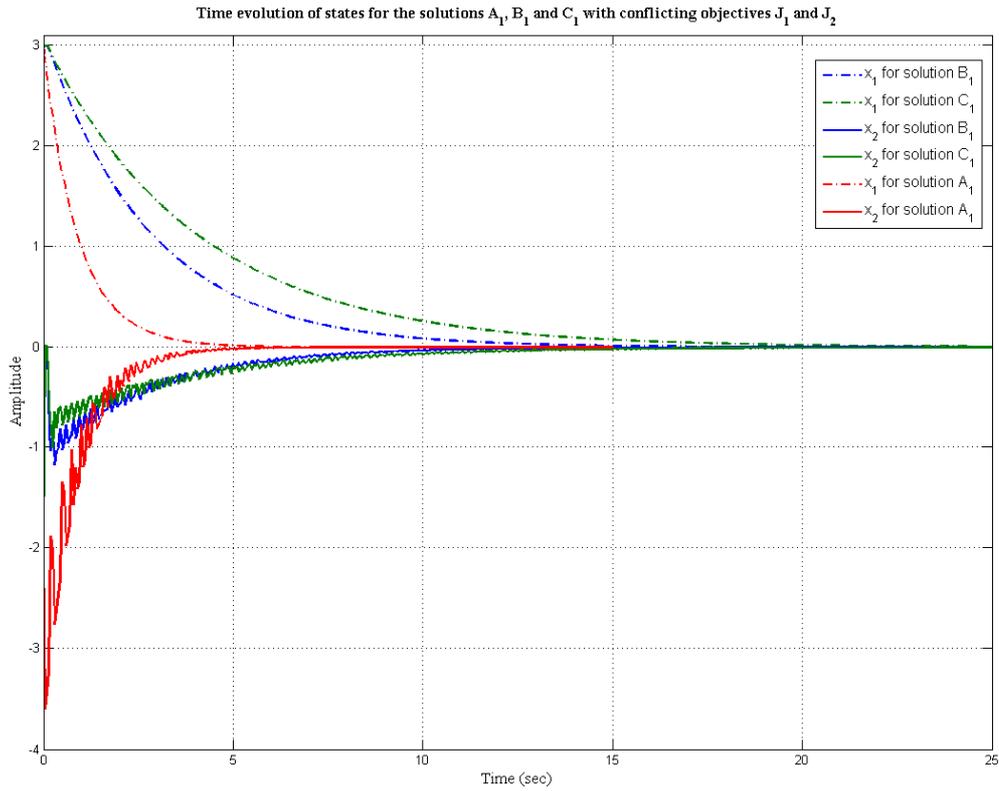

Fig. 5. Time evolution of states for solutions $A_1$, $B_1$ and $C_1$ with the conflicting objectives $J_1$ and $J_2$

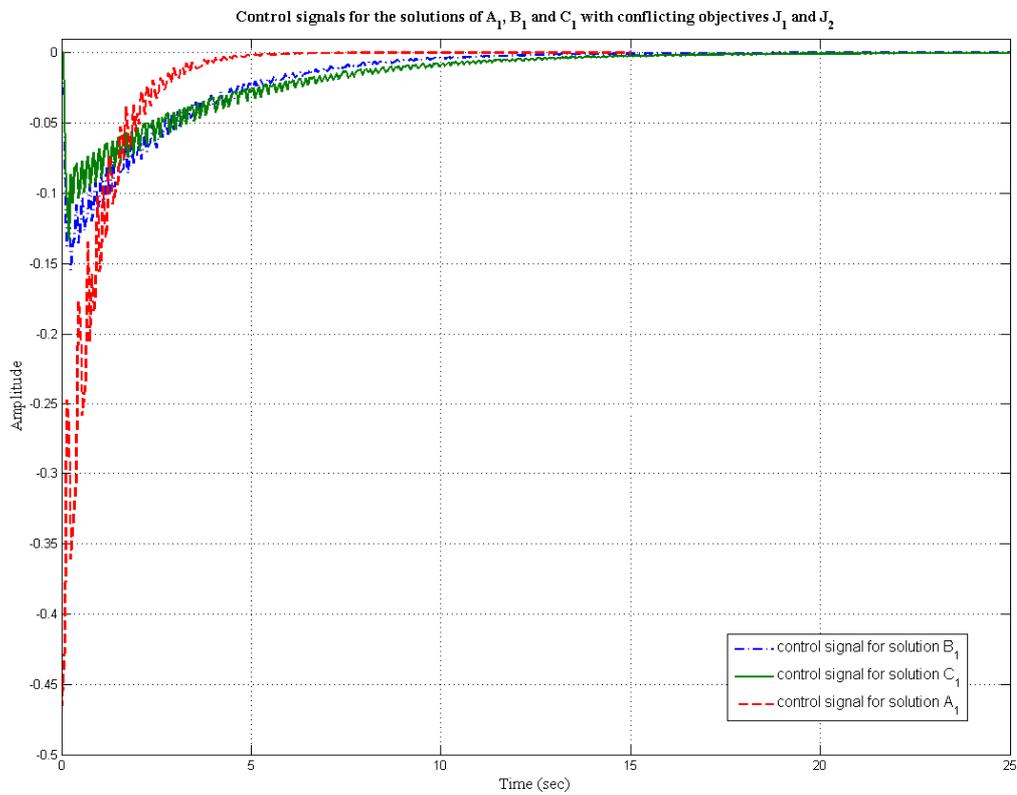

Fig. 6. Plot of control signals for the solutions $A_1$, $B_1$ and $C_1$.

Similarly, the Pareto-front for the objective functions $J_3$ and $J_2$ has been shown in Fig. 7. Similar to the previous case, two extreme solutions and one intermediate solution is taken and their controller values and time domain performances are shown in Table 2 and Fig. 8-9 respectively. From Fig. 8 and 9 it is clear that making state transition smooth, lead to oscillatory control signals, because the control signal changes rapidly so that the states can track a smooth curve. But smoothing of states cause them to take a long time to settle down as the settling time is not penalized in the designing technique.

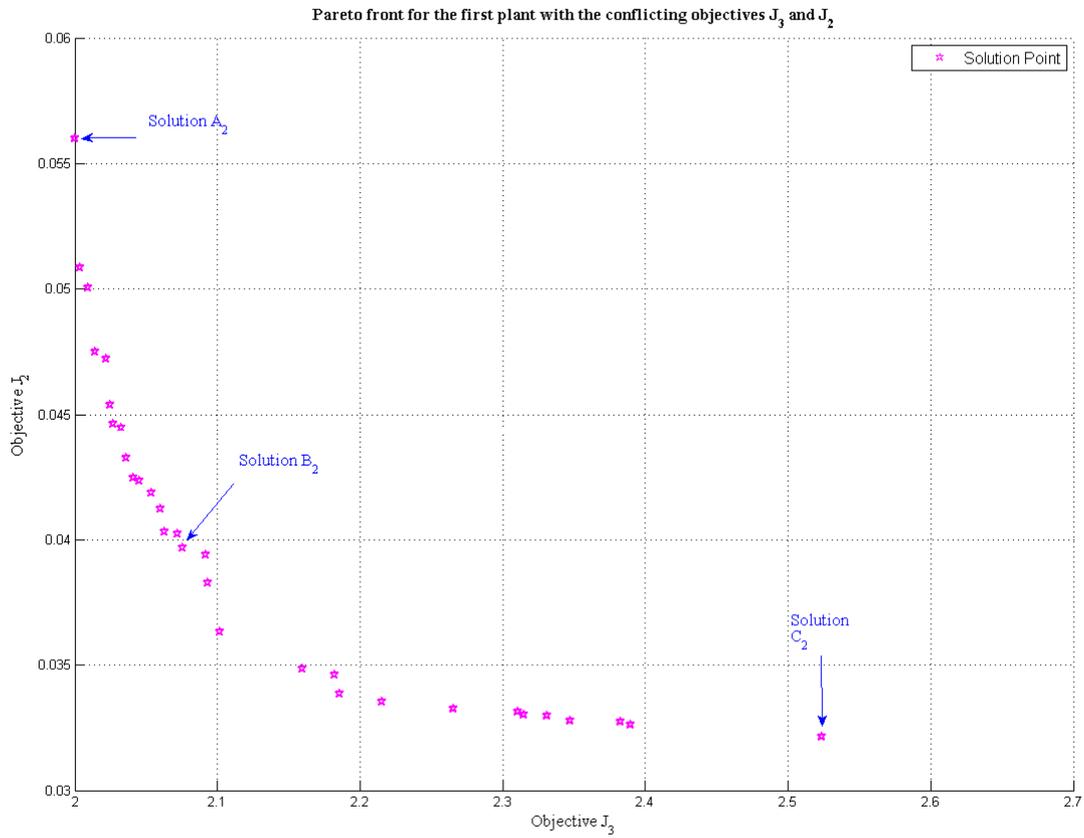

Fig. 7. Pareto optimal front for the conflicting objective functions $J_3$ and $J_2$

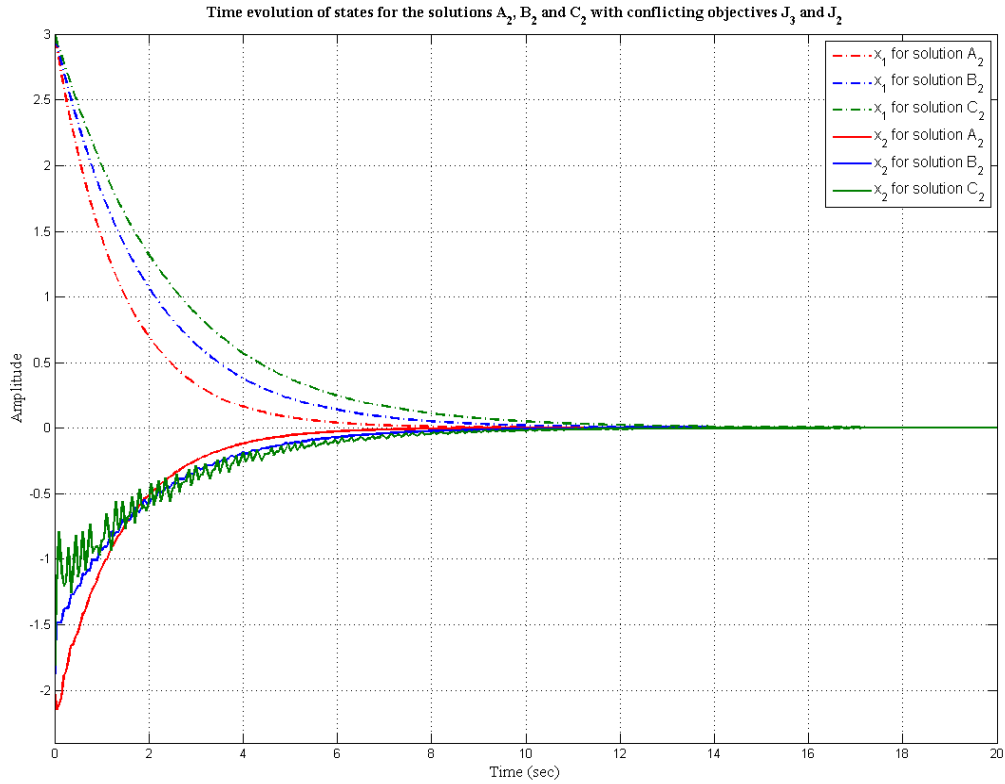

Fig. 8. Time evolution of states for solutions $A_2$, $B_2$ and $C_2$ with the conflicting objectives $J_3$ and $J_2$

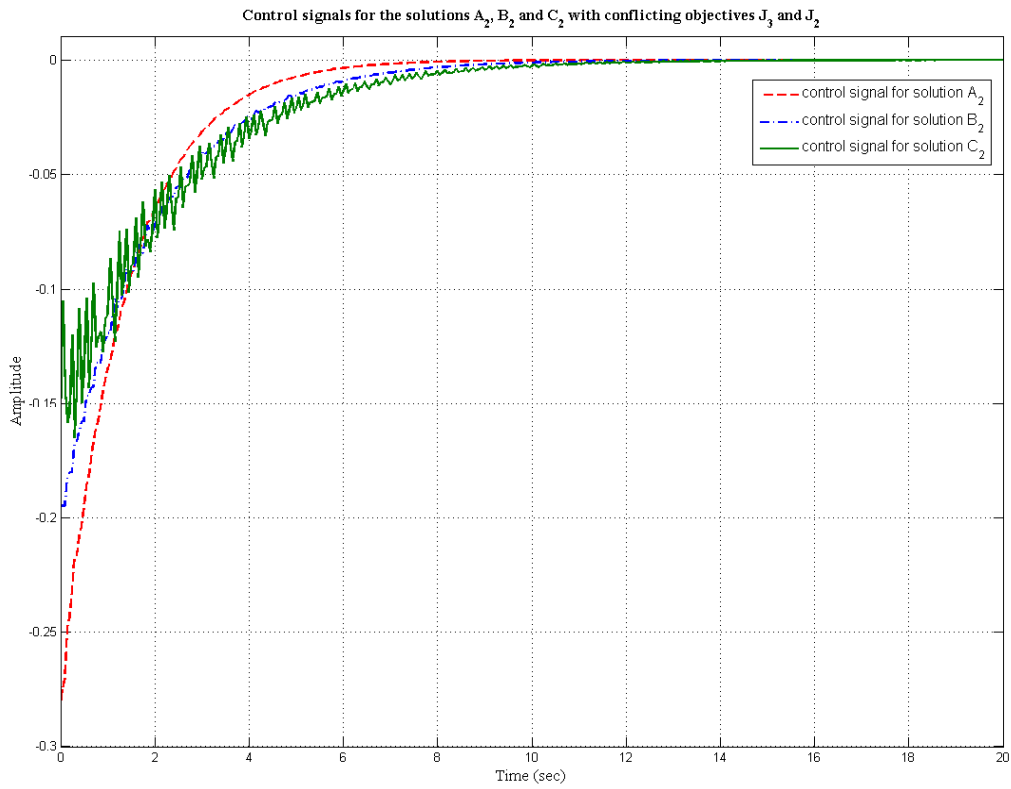

Fig. 9. Plot of control signals for the solutions $A_2$, $B_2$ and $C_2$

Now for the third type of design technique (as discussed in Section 3.3) the Pareto plot obtained after optimization is shown in Fig. 10 and three solution points $A_3$, $B_3$ and $C_3$ are taken for further investigation.

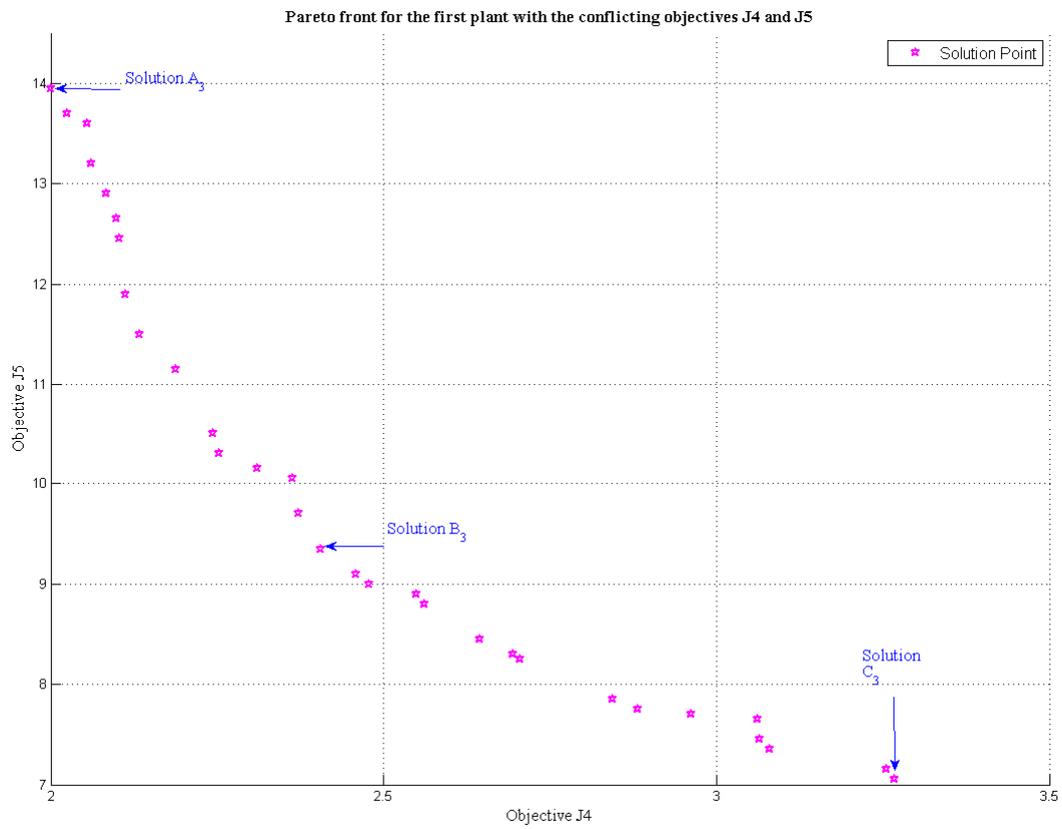

Fig. 10. Pareto optimal front for the conflicting objective functions $J_4$ and $J_5$

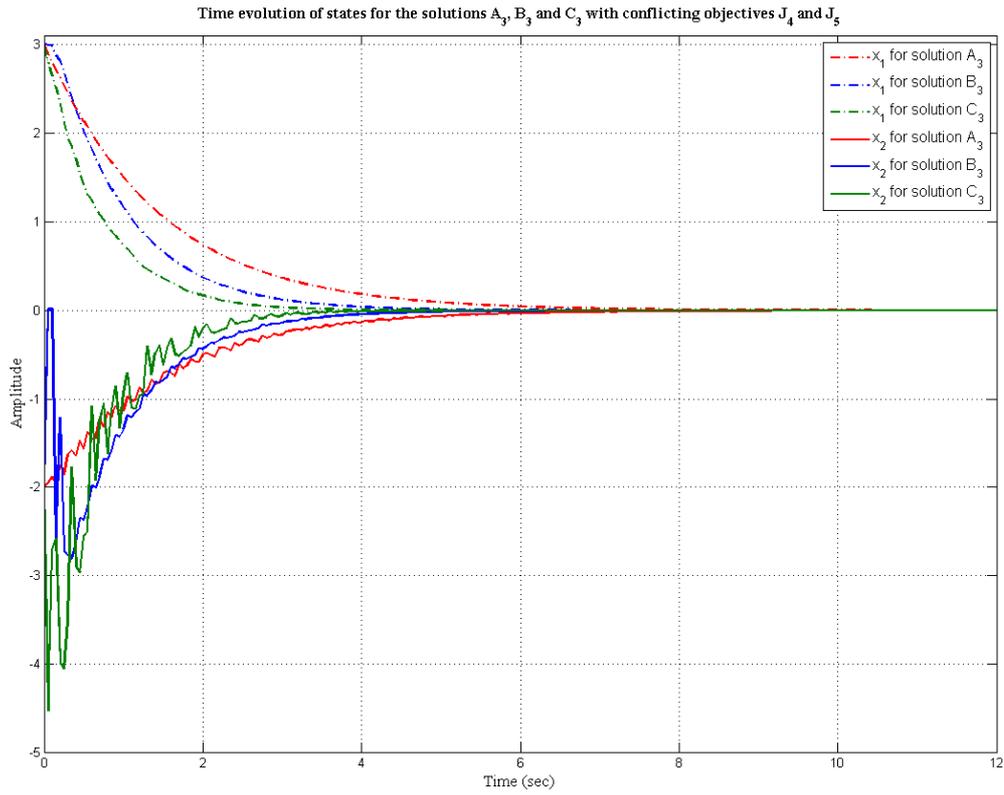

Fig. 11. Time evolution of states for solutions $A_3$, $B_3$ and $C_3$ with the conflicting objectives $J_4$ and $J_5$

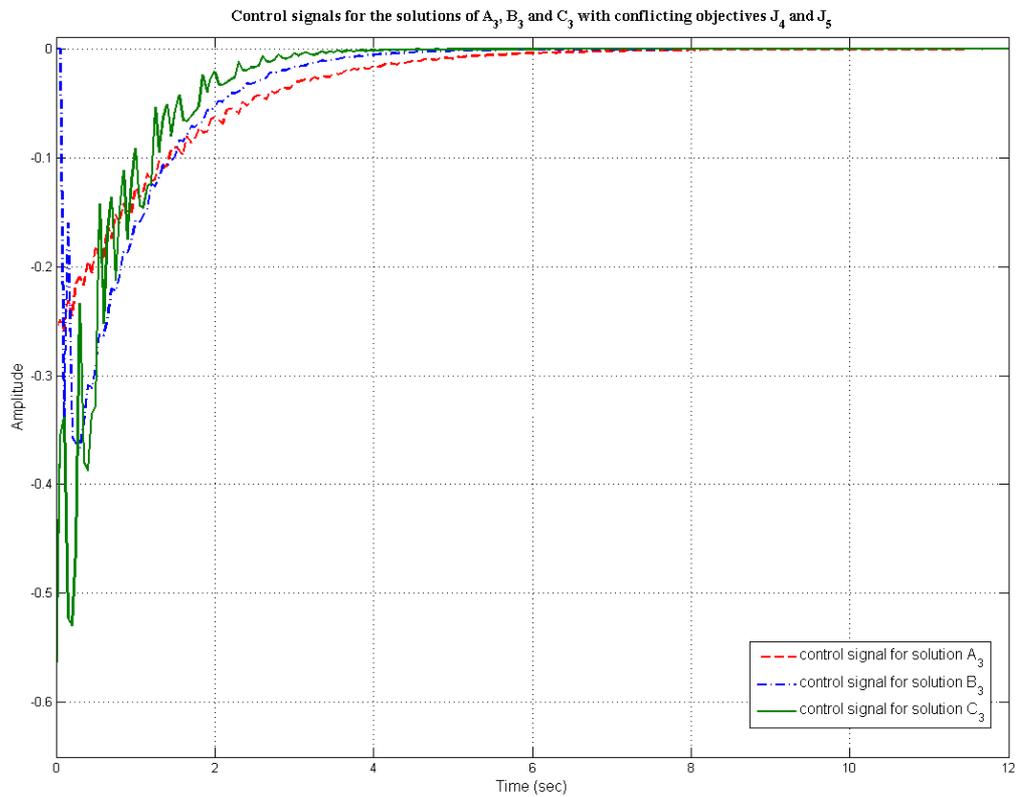

Fig. 12. Plot of control signals for the solutions $A_3$, $B_3$ and $C_3$

Figs. 11 and 12 show the time response and the control signal of the representative solutions respectively. It is evident from the Table 2 and also Fig. 11 that response for $C_3$ is quicker than others and response of $A_3$ takes longest time to settle. But states for $C_3$ attain a higher peak value than the states for $A_3$. Fig. 12 shows that more control signal is needed for $C_3$ than for $A_3$. Hence it is proved from the graphs that for quicker response overshoot increases, for which more control signal is needed to settle down in a short time. Hence overshoot and settling time are two conflicting objectives that need to be traded-off for an optimal controller design.

## 4.2. Double Integrator plant

The unstable system, considered next has two poles at the origin of complex s-plane which represents a double integrator process [35].

$$\dot{x}(t) = \begin{bmatrix} 0 & 1 \\ 0 & 0 \end{bmatrix} x(t) + \begin{bmatrix} 0 \\ 1 \end{bmatrix} u(t) \tag{18}$$

By taking sampling period as 0.01s, we can get discrete time plant having the system matrices given by (19)

$$F = \begin{bmatrix} 1 & 0.01 \\ 0 & 1 \end{bmatrix}, \ G = \begin{bmatrix} 0.0001 \\ 0.01 \end{bmatrix} \tag{19}$$

Multi-objective optimization for the predictive controller design with different objective functions yields controller parameters given in Table 3 while considering the same initial condition as before.

Table 3:

Multi-objective optimization results amongst the objective functions for the double integrator plant

| Design trade-off amongst objective functions | Solution points | Objective functions | | | | | State feedback controller gains | | | | | |
|---|---|---|---|---|---|---|---|---|---|---|---|---|
| | | $J_1$ | $J_2$ | $J_3$ | $J_4$ | $J_5$ | $K_{11}$ | $K_{12}$ | $K_{21}$ | $K_{22}$ | $K_{31}$ | $K_{32}$ |
| $J_1$ and $J_2$ | $A_1$ | 4.716 | 5.623 | - | - | - | -1.534 | -1.650 | -1.179 | -2.644 | -3.497 | -2.730 |
| | $B_1$ | 6.525 | 1.251 | - | - | - | -0.994 | -1.670 | -0.658 | -1.279 | -1.374 | -2.164 |
| | $C_1$ | 20.022 | 1.094 | - | - | - | -0.253 | -0.741 | -0.439 | -0.970 | -0.139 | -0.719 |
| $J_3$ and $J_2$ | $A_2$ | - | 8.816 | 1.939 | - | - | -2.660 | -1.593 | -0.814 | -1.233 | -2.401 | -1.476 |
| | $B_2$ | - | 2.510 | 1.946 | - | - | -1.678 | -1.624 | -0.917 | -1.480 | -1.540 | -1.721 |
| | $C_2$ | - | 1.076 | 1.954 | - | - | -0.309 | -0.686 | -0.417 | -0.962 | -0.161 | -0.627 |
| $J_4$ and $J_5$ | $A_3$ | - | - | - | 2.000 | 5.900 | -0.949 | -1.425 | -0.598 | -0.867 | -0.707 | -1.366 |
| | $B_3$ | - | - | - | 2.024 | 5.280 | -1.086 | -1.644 | -0.355 | -1.301 | -2.343 | -1.723 |
| | $C_3$ | - | - | - | 2.176 | 4.570 | -1.784 | -1.276 | -1.127 | -2.745 | -3.318 | -2.731 |

Fig. 13 shows the Pareto front for the trade-off between the two objectives $J_1$ and $J_2$ for the double integrator plant. Fig. 14 and 15 shows the time evolution of the states and the control signals respectively for the three representative solutions as labelled on the Pareto front. Similar to the dc motor plant in the previous case, the solution $A_1$ has the fastest settling time but requires a much larger control signal. Solution $C_1$ has the longest settling time among the three solutions but requires much lower value of the control signal. The key difference between this example and the previous one of the dc motor is that there are a lot of oscillations in the control signal. This is due to the nature of the double integrator plant which is inherently unstable. Similar example of unstable system stabilization with fractional order and integer order PID type controllers over network as in Pan *et al.* [36] have shown that the control signal may be jittery even though the time response curves are smooth. This is due the fact that the controller produces additional control signals after the packet drop occurs and to compensate for this missing measurement or control values the manipulated variable is perturbed violently. This phenomenon is intrinsically different from that due to measurement noise and cannot be removed by conventional derivative filtering or similar techniques. It can be seen that for solution $A_1$ although there is an improvement in the settling time, there are much larger oscillations in the control signal. So, if this is a mechanical system, then it is not really advisable to go with solution $A_1$ as large frequent variations in control signal would result in dithering of the actuator which is undesirable.

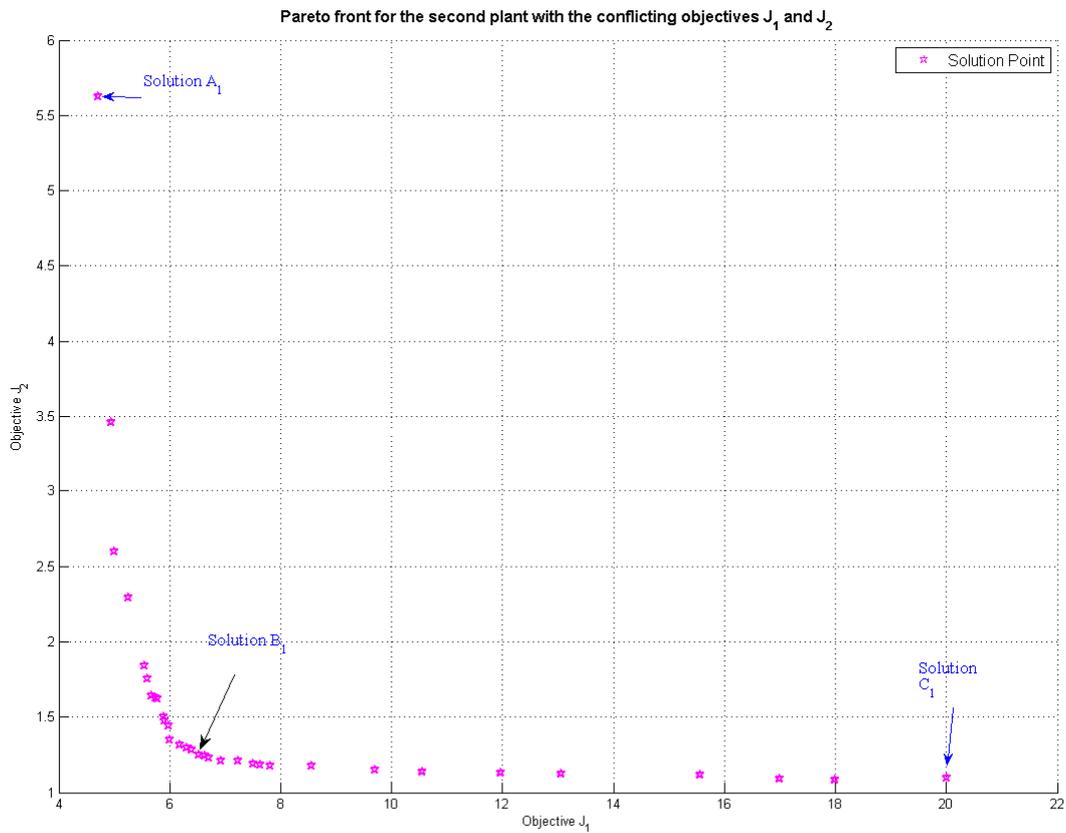

Fig. 13. Pareto optimal front for the conflicting objective functions $J_1$ and $J_2$

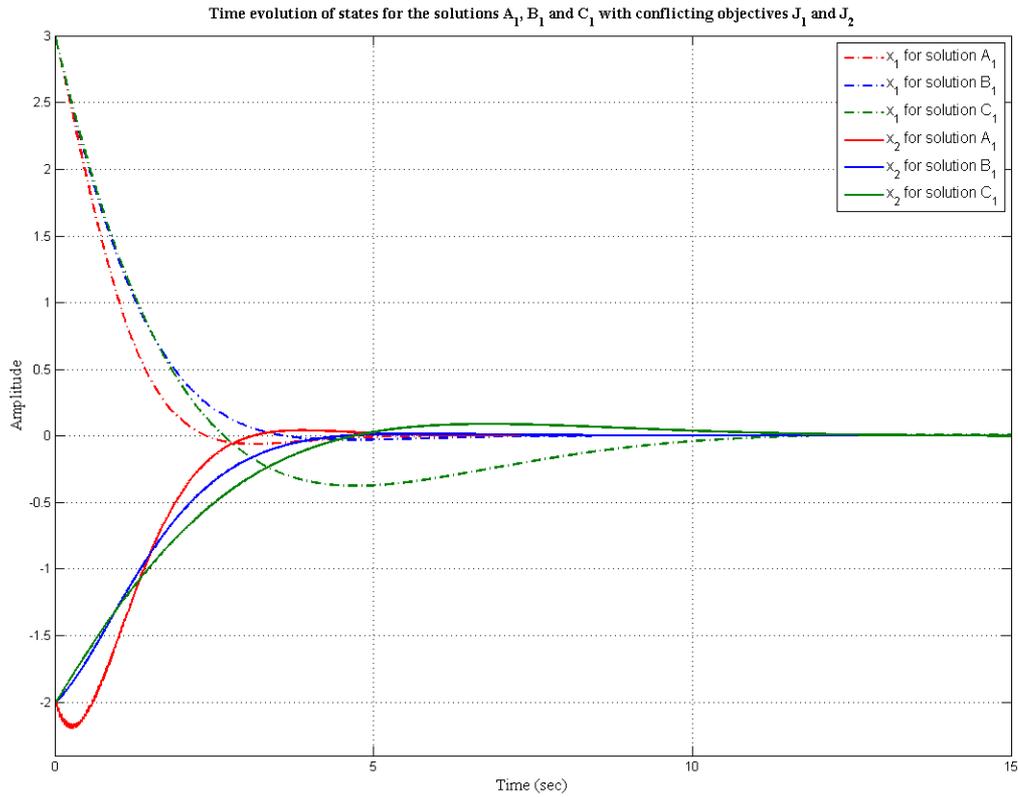

Fig. 14. Time evolution of states for solutions $A_1$, $B_1$ and $C_1$ with the conflicting objectives $J_1$ and $J_2$

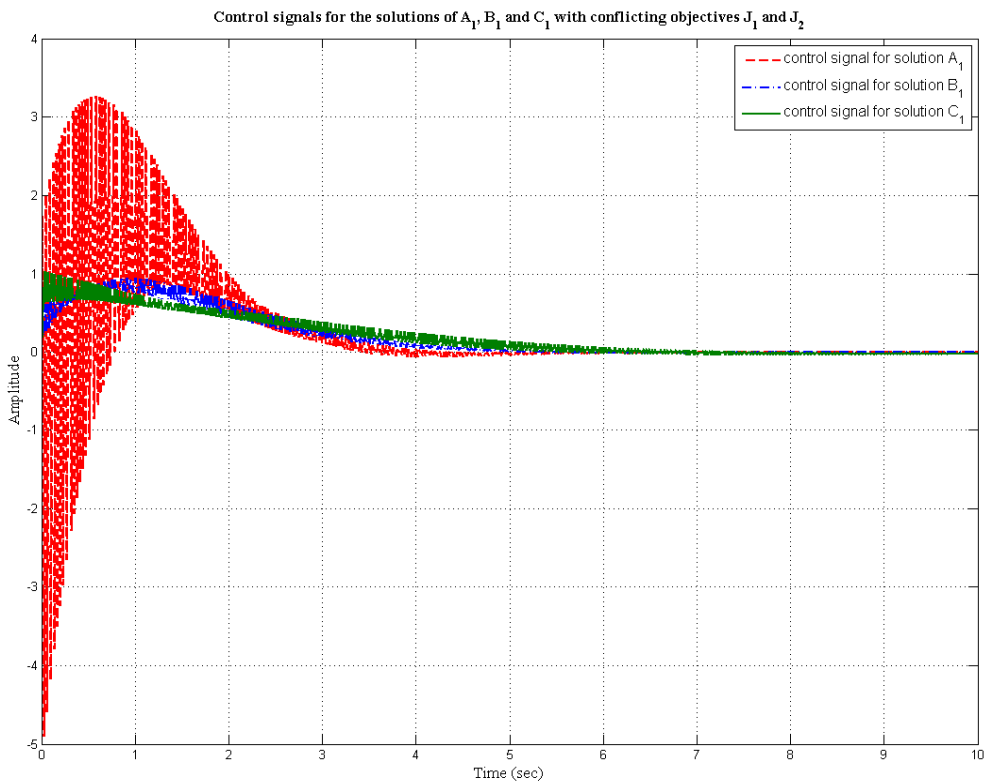

Fig. 15. Plot of control signals for the solutions $A_1$, $B_1$ and $C_1$.

Fig. 16 shows the Pareto front for the trade-off between $J_2$ and $J_3$ with three representative solutions labelled on the front. Figs. 17 and 18 show the time response and the control signals of these representative solutions respectively. From Fig. 18, it is found that the control signal for $A_2$ is the highest with very high oscillations in the signal as compared to solutions $B_2$ and $C_2$. This infact is also corroborated by the position of $A_2$ on the Pareto front in Fig. 16. The state transitions of all the solutions are mostly smooth in Fig. 17 and the minute differences is almost negligible at a glance. This is also due to the fact that the range of x-axis (objective $J_3$) in Fig. 16 is much smaller than the range of the y-axis (objective $J_2$). However as a consequence of the larger control signal in $A_2$, the solution has automatically been the fastest to settle. But since there is no constraint on the peak overshoot, the same solution has the maximum peak overshoot as compared to the other ones.

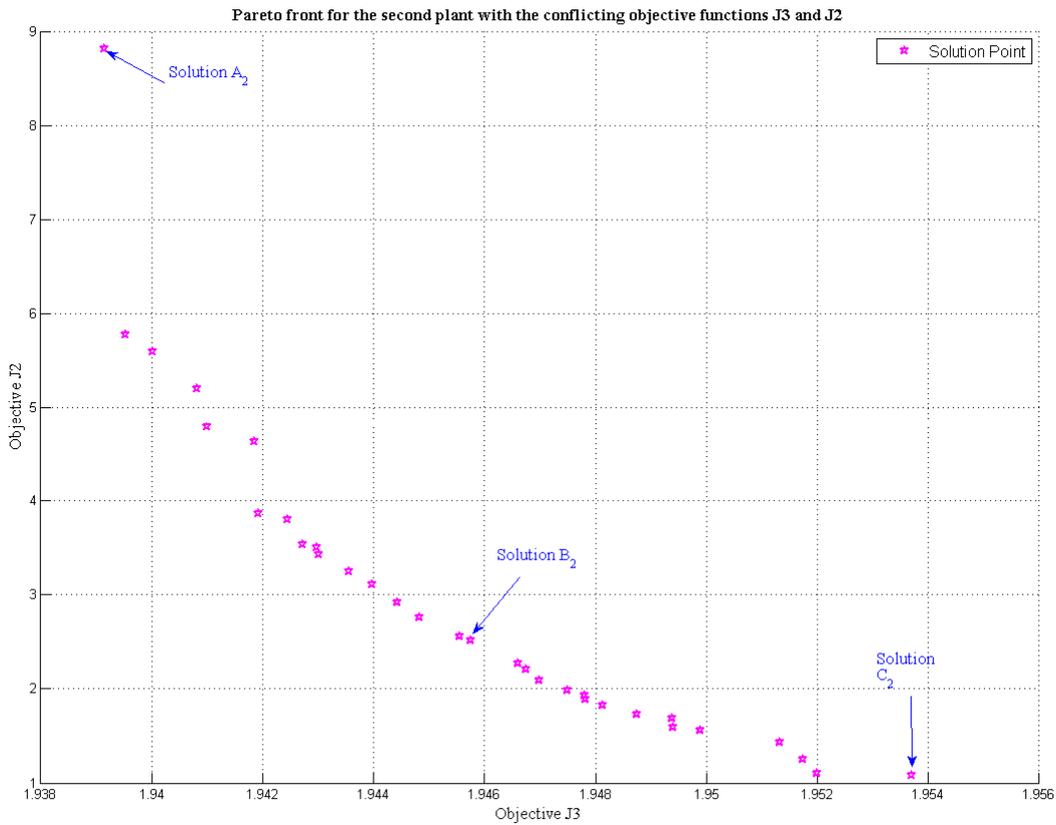

Fig. 16. Pareto optimal front for the conflicting objective functions $J_3$ and $J_2$

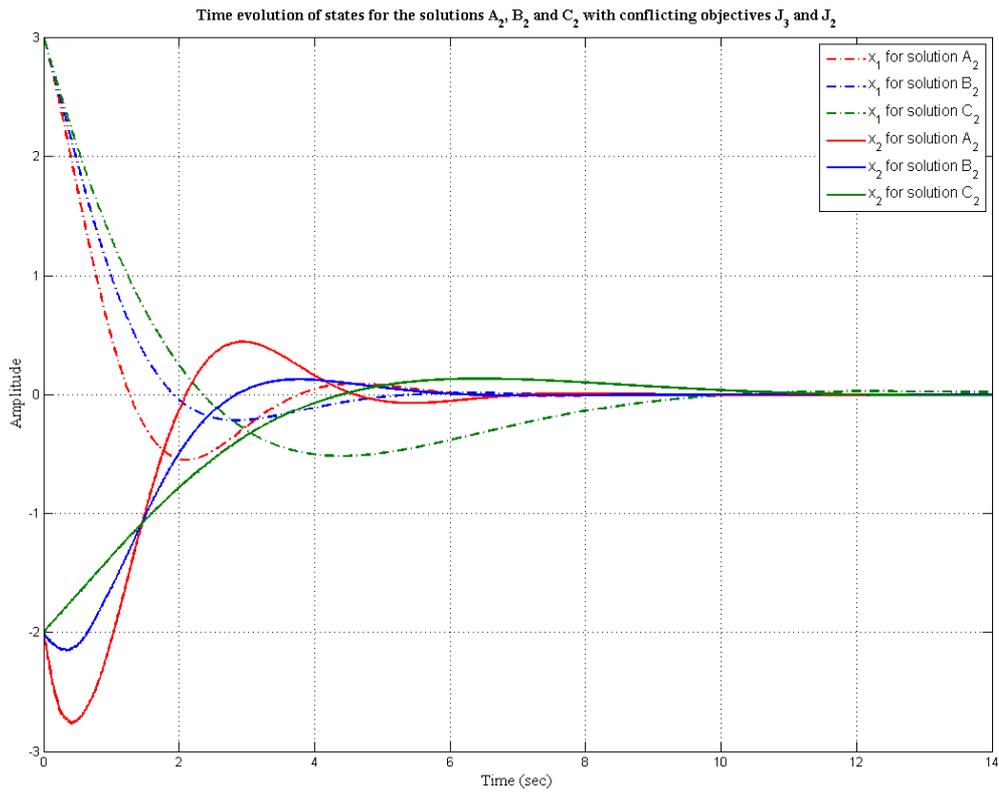

Fig. 17. Time evolution of states for solutions $A_2$, $B_2$ and $C_2$ with the conflicting objectives $J_3$ and $J_2$

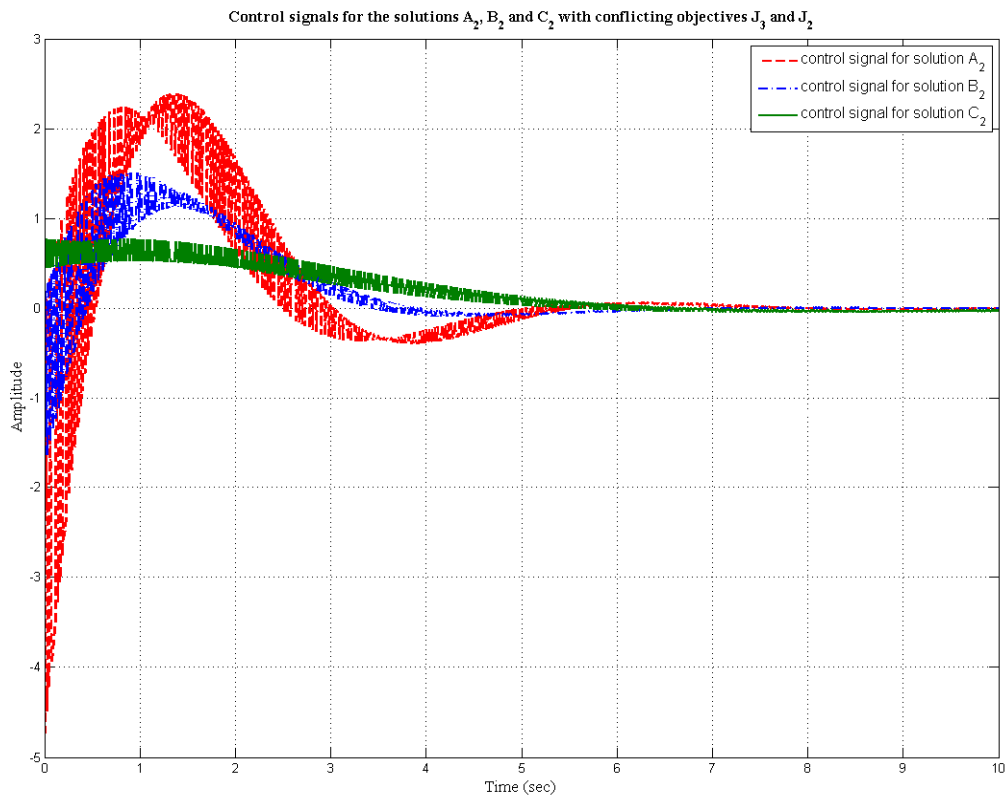

Fig. 18. Plot of control signals for the solutions $A_2$, $B_2$ and $C_2$

Fig. 19 shows the Pareto trade-off for the two conflicting objective functions $J_4$ and $J_5$ for the double integrator plant. Figs. 20 and 21 show the plots of the time response and the control signal respectively of the representative solutions as labelled on the Pareto front. As observed from Fig. 20, the peak overshoot is higher for solution $C_3$ as compared to $A_3$ but the settling time of $C_3$ is faster than that of $A_3$. Since the x-axis (objective $J_4$) of the Pareto front in Fig. 19 has a very small range, hence the difference in the settling time of the two solutions $C_3$ and $A_3$ are very small. Now a higher peak overshoot or faster settling time implies a higher control signal. This can be verified from Fig. 21 where the control signal of solution $C_3$ is the highest.

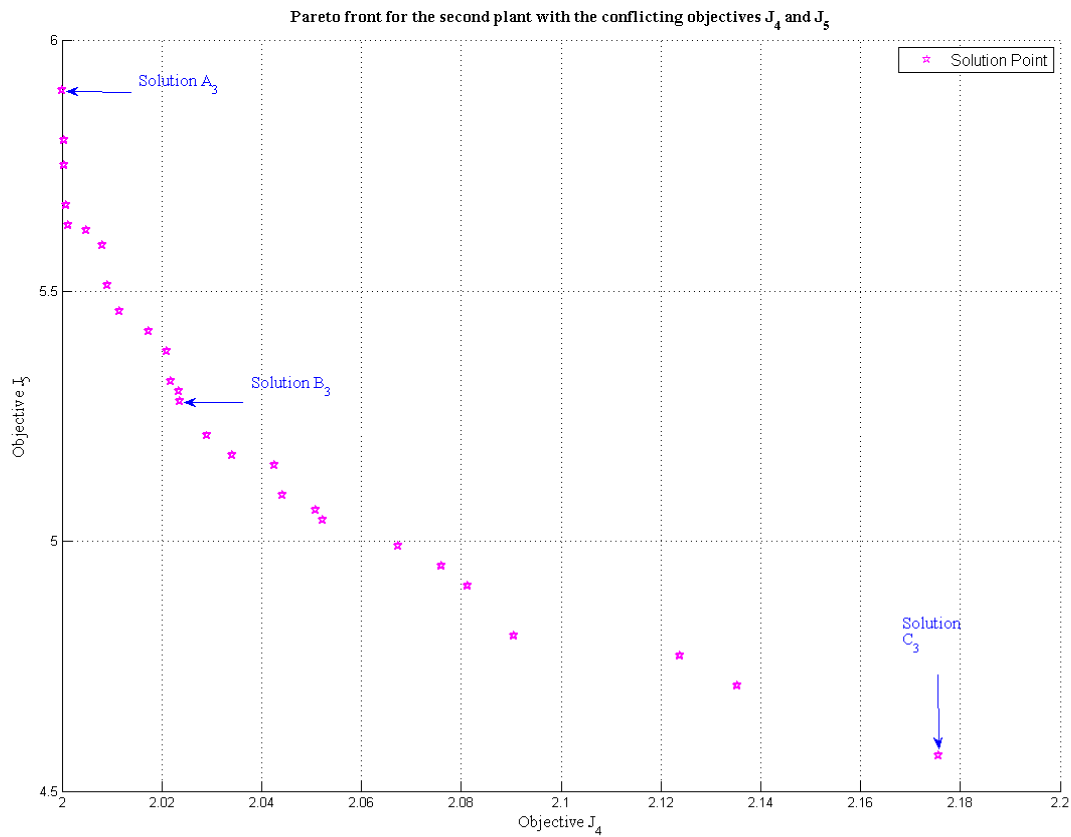

Fig. 19. Pareto optimal front for the conflicting objective functions $J_4$ and $J_5$

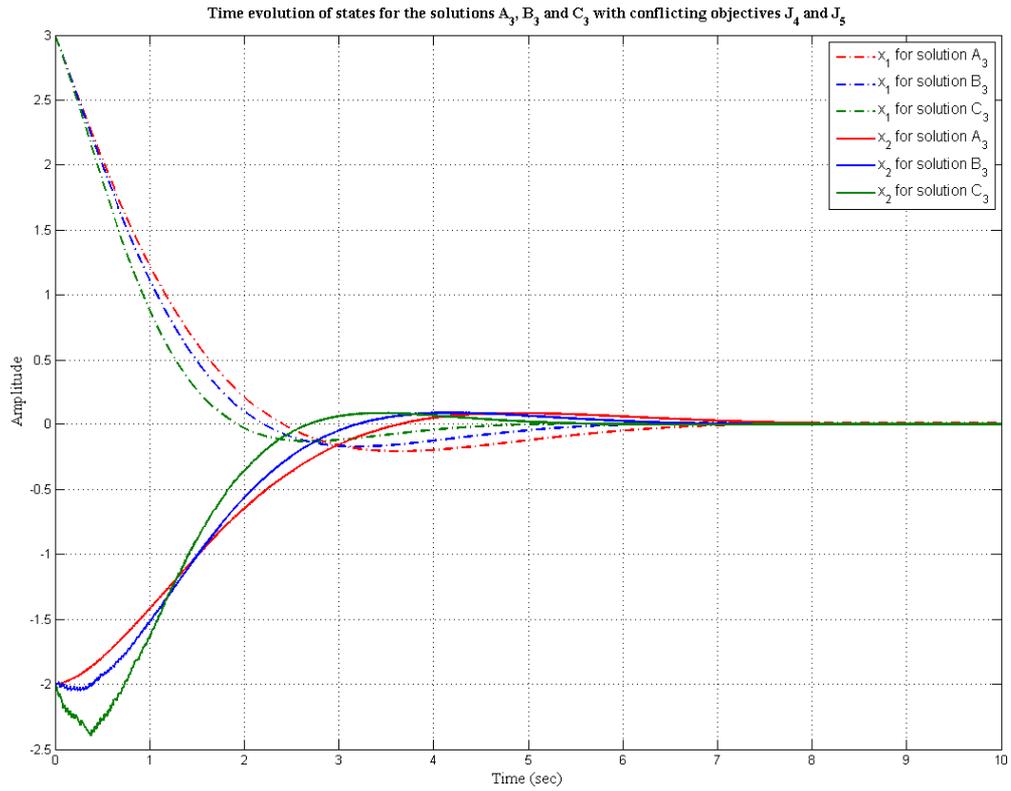

Fig. 20. Time evolution of states for solutions $A_3$, $B_3$ and $C_3$ with the conflicting objectives $J_4$ and $J_5$

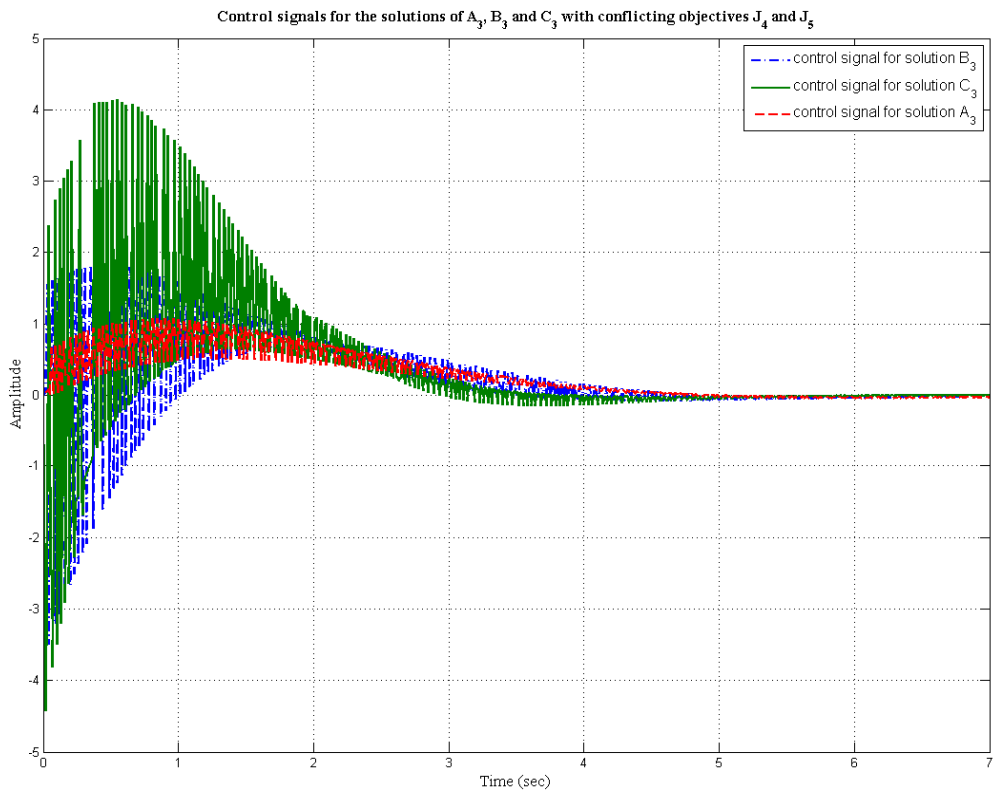

Fig. 21. Plot of control signals for the solutions $A_3$, $B_3$ and $C_3$.

### 4.3. Inverted pendulum plant

The continuous time plant, considered here has one stable pole and one unstable pole from which are at equal distance from the origin. It has been shown in [37] that such a simplified model represents an inverted pendulum process.

$$\dot{x}(t) = \begin{bmatrix} 0 & 1 \\ 1 & 0 \end{bmatrix} x(t) + \begin{bmatrix} 0 \\ 1 \end{bmatrix} u(t) \qquad (20)$$

By choosing the sampling time as 0.05s, we can get discrete time plant having the system matrices given by (21)

$$F = \begin{bmatrix} 1.0013 & 0.05 \\ 0.05 & 1.0013 \end{bmatrix}, G = \begin{bmatrix} 0.0013 \\ 0.05 \end{bmatrix} \qquad (21)$$

For the above system, multi-objective optimization with LMI criteria yields stabilizing Pareto optimal gains as reported in Table 4 while taking different objective functions for simulation. Fig. 22 shows the Pareto front for the inverted pendulum for the two contradictory objectives $J_1$ and $J_2$. Fig. 23 and 24 show the state responses and the control signal respectively for the inverted pendulum for the three representative cases as labelled in the Pareto front. Similar results can be observed as in the previous two cases with solution $A_1$ having the fastest settling time and solution $C_1$ having the slowest. Solution $B_1$ has a settling time in between the two extremes. Since this is an unstable plant, the control signal for $A_1$ has significantly more oscillations than solutions $B_1$ and $C_1$. But nevertheless, these oscillations are much lower than that in the double integrator case.

Table 4:

Multi-objective optimization results amongst the objective functions for the inverted pendulum plant

| Design trade-off amongst objective functions | Solution points | Objective functions | | | | | State feedback controller gains | | | | | |
|---|---|---|---|---|---|---|---|---|---|---|---|---|
| | | $J_1$ | $J_2$ | $J_3$ | $J_4$ | $J_5$ | $K_{11}$ | $K_{12}$ | $K_{21}$ | $K_{22}$ | $K_{31}$ | $K_{32}$ |
| $J_1$ and $J_2$ | $A_1$ | 5.579 | 3.597 | - | - | - | -2.486 | -2.322 | -2.068 | -1.047 | -1.847 | -1.878 |
| | $B_1$ | 6.389 | 2.354 | - | - | - | -2.370 | -2.229 | -1.468 | -1.028 | -1.929 | -1.843 |
| | $C_1$ | 8.942 | 2.002 | - | - | - | -2.010 | -2.048 | -1.288 | -1.052 | -1.797 | -1.804 |
| $J_3$ and $J_2$ | $A_2$ | - | 26.993 | 1.938 | - | - | -3.115 | -1.015 | -1.914 | -0.872 | -2.749 | -1.198 |
| | $B_2$ | - | 5.676 | 1.964 | - | - | -2.644 | -1.746 | -1.927 | -1.280 | -2.201 | -1.686 |
| | $C_2$ | - | 2.002 | 2.000 | - | - | -1.976 | -2.021 | -1.911 | -1.945 | -1.928 | -1.971 |
| $J_4$ and $J_5$ | $A_3$ | - | - | - | 2.000 | 6.300 | -2.982 | -1.571 | -2.511 | -2.130 | -1.859 | -3.151 |
| | $B_3$ | - | - | - | 2.028 | 6.050 | -2.906 | -1.569 | -2.492 | -2.142 | -1.883 | -3.161 |
| | $C_3$ | - | - | - | 2.251 | 5.650 | -3.005 | -1.574 | -2.552 | -2.176 | -1.867 | -3.158 |

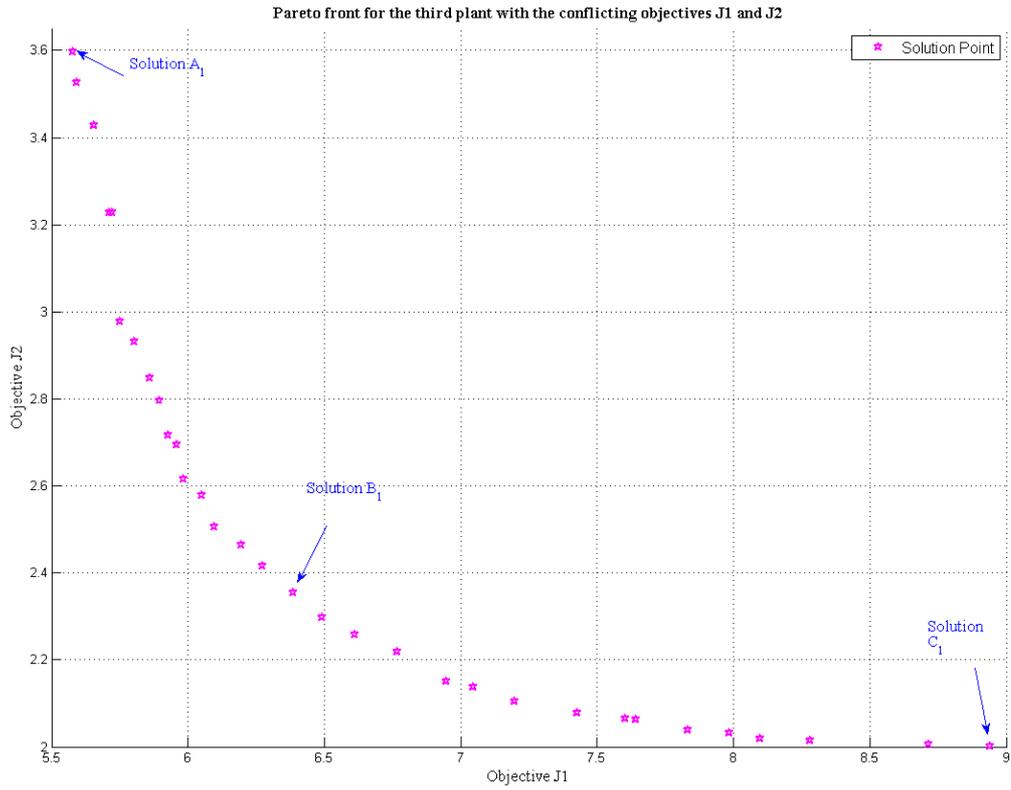

Fig. 22. Pareto optimal front for the conflicting objective functions $J_1$ and $J_2$

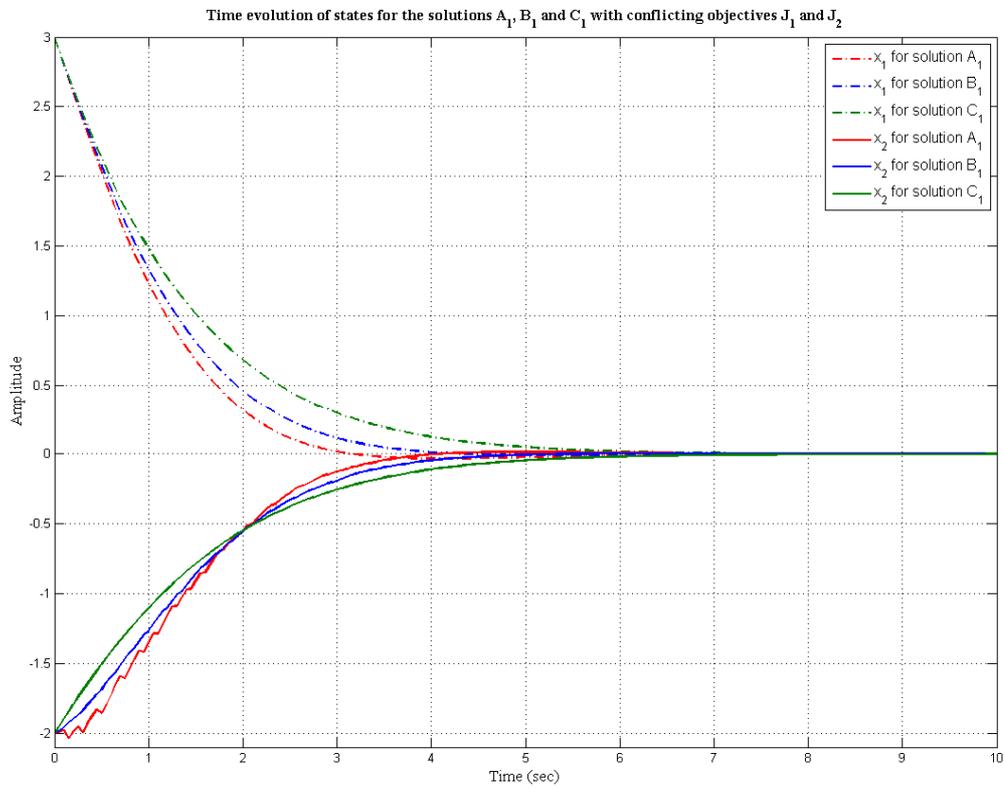

Fig. 23. Time evolution of states for solutions $A_1$, $B_1$ and $C_1$ with the conflicting objectives $J_1$ and $J_2$

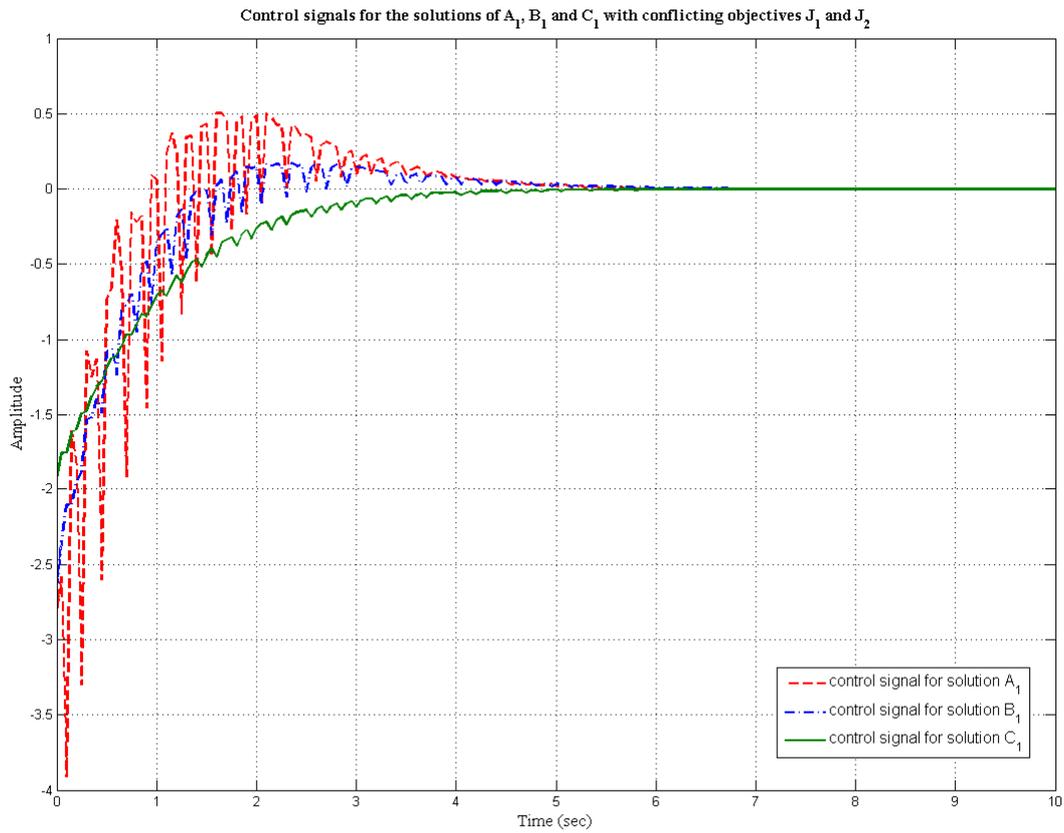

Fig. 24. Plot of control signals for the solutions $A_1$, $B_1$ and $C_1$.

Fig. 25 shows the Pareto front for the two contradictory objectives $J_2$ and $J_3$ for the inverted pendulum plant with three representative solutions labelled. Figs. 26 and 27 show the time response plots of these solutions for these representative cases. Similar results are seen as in the previous two plants (the dc motor and the double integrator). The increase in control signal for $A_2$ automatically makes it the fastest to settle and gives higher overshoot than the others. The oscillations in the control signal for the inverted pendulum are less than that of the double integrator plant.

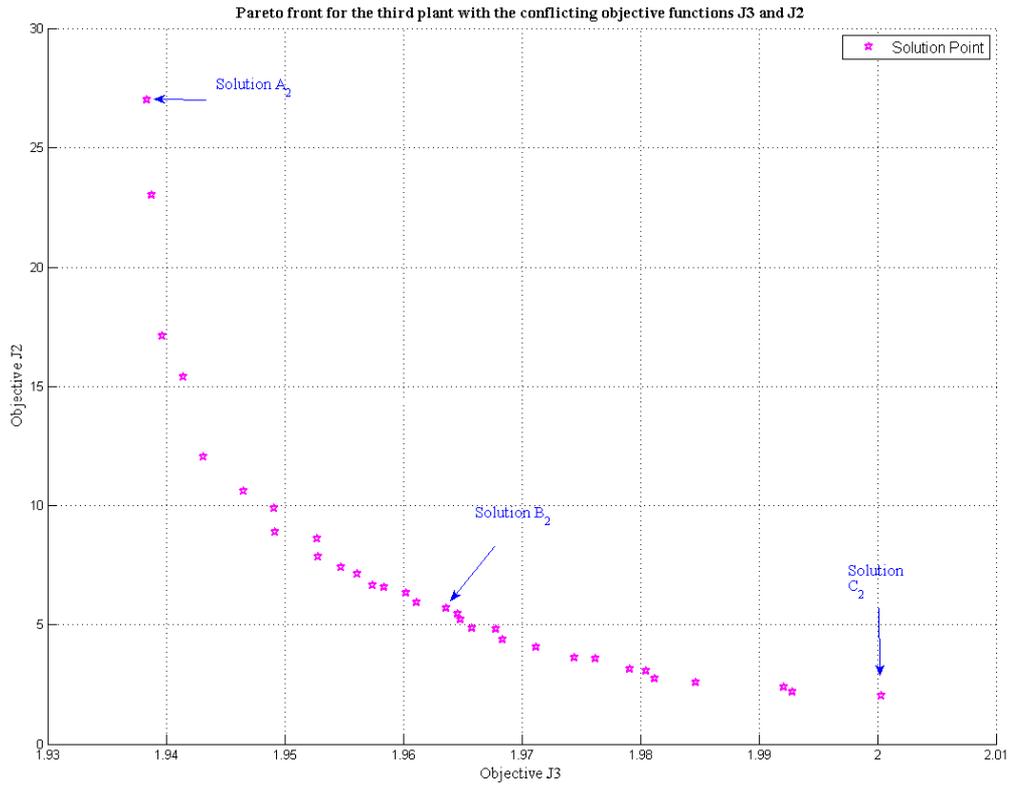

Fig. 25. Pareto optimal front for the conflicting objective functions $J_3$ and $J_2$

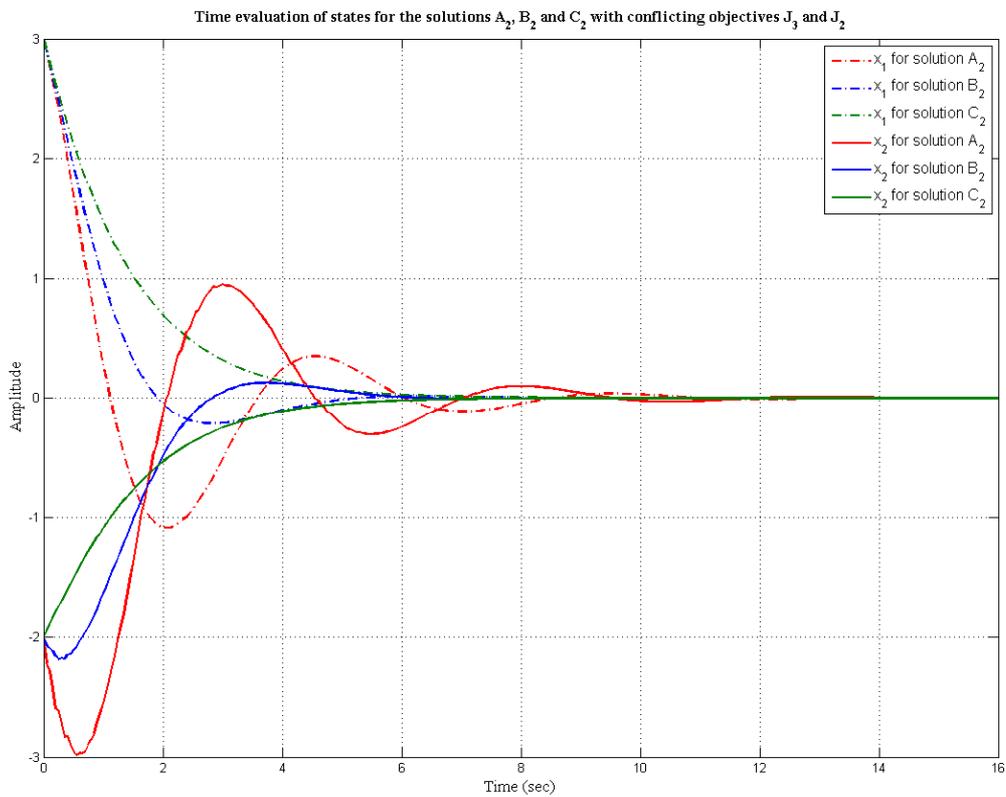

Fig. 26. Time evolution of states for solutions $A_2$, $B_2$ and $C_2$ with the conflicting objectives $J_3$ and $J_2$

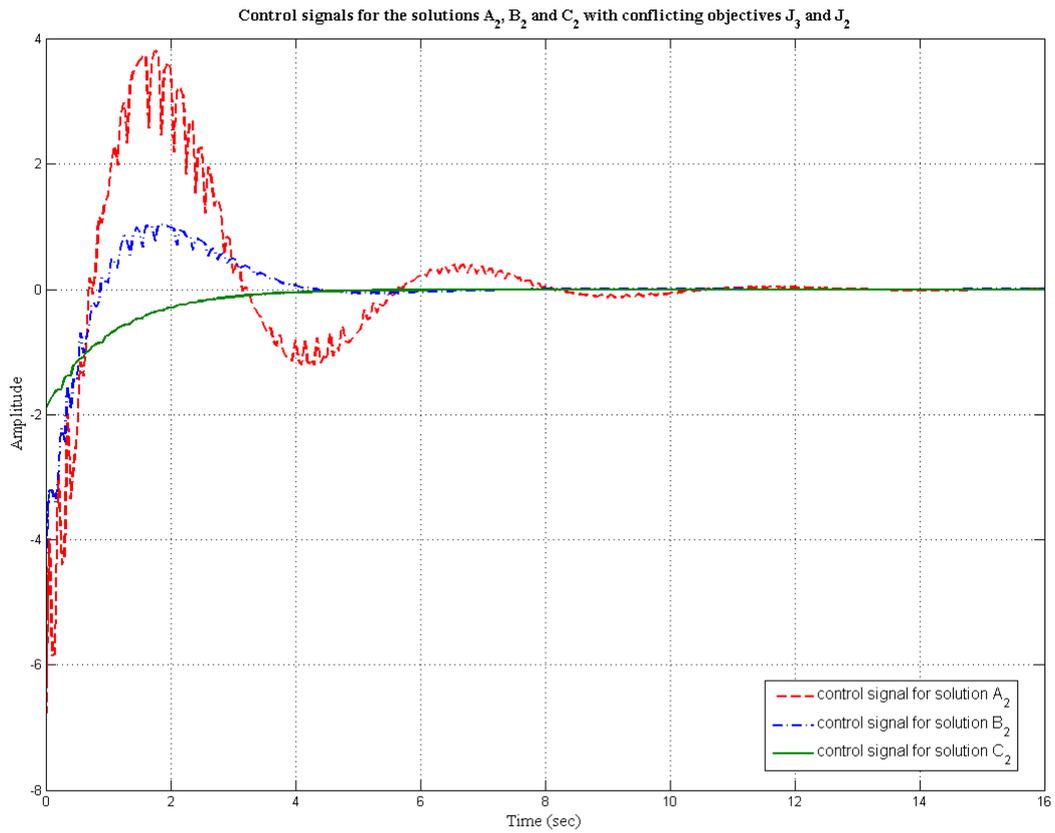

Fig. 27. Plot of control signals for the solutions $A_2$, $B_2$ and $C_2$

Fig. 28 shows the Pareto front with representative solutions for the conflicting objective functions $J_4$ and $J_5$ for the inverted pendulum plant. Figs. 29 and 30 show the time domain evolution of the states and the control signals respectively for the representative solutions. Similar characteristics of the solutions can be seen as in the previous two cases of the dc motor and the double integrator plant. The control signals are jittery for all the obtained solutions in this case.

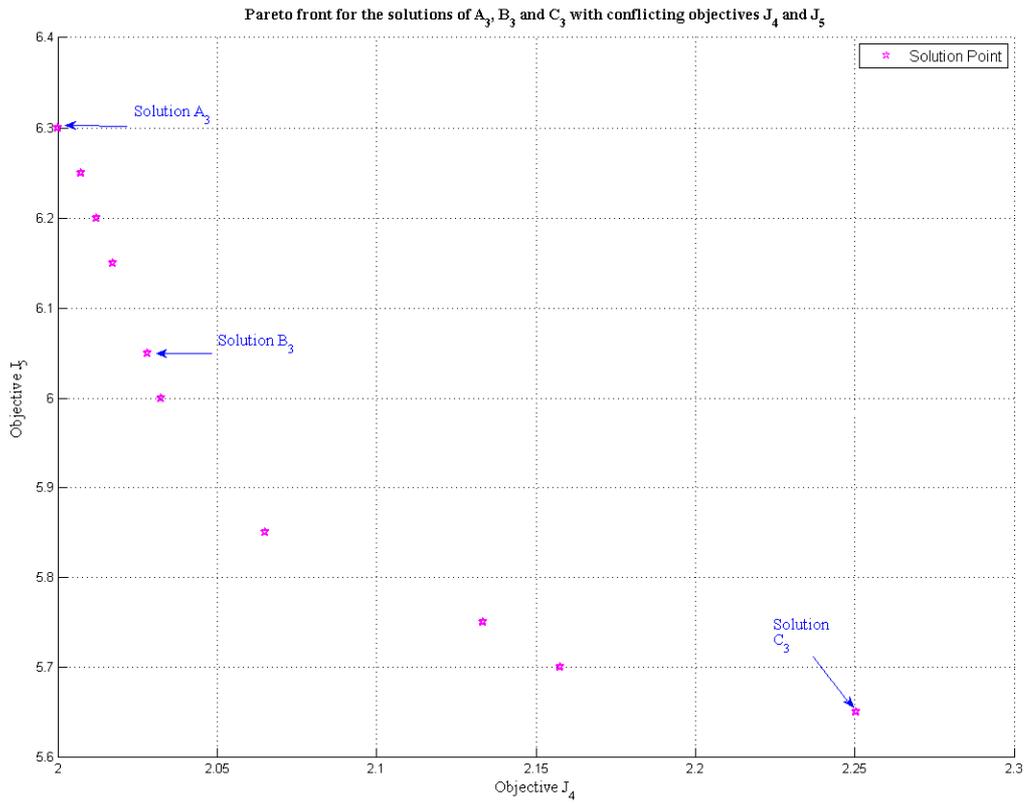

Fig. 28. Pareto optimal front for the conflicting objective functions $J_4$ and $J_5$

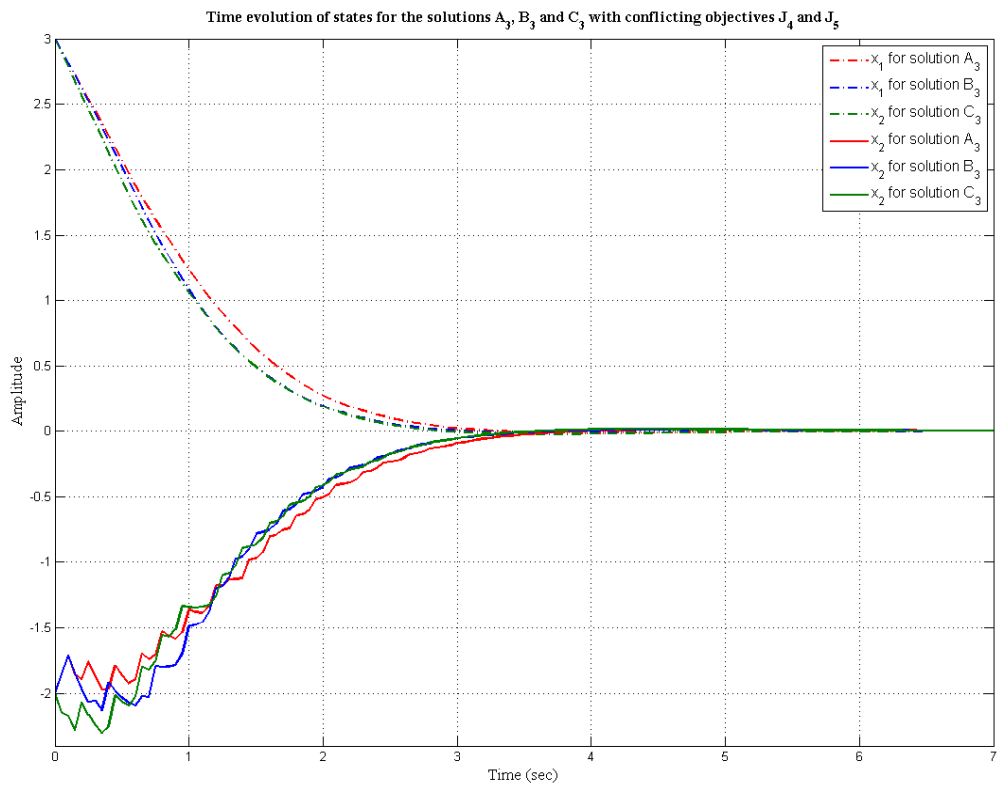

Fig. 29. Time evolution of states for solutions $A_3$, $B_3$ and $C_3$ with the conflicting objectives $J_4$ and $J_5$

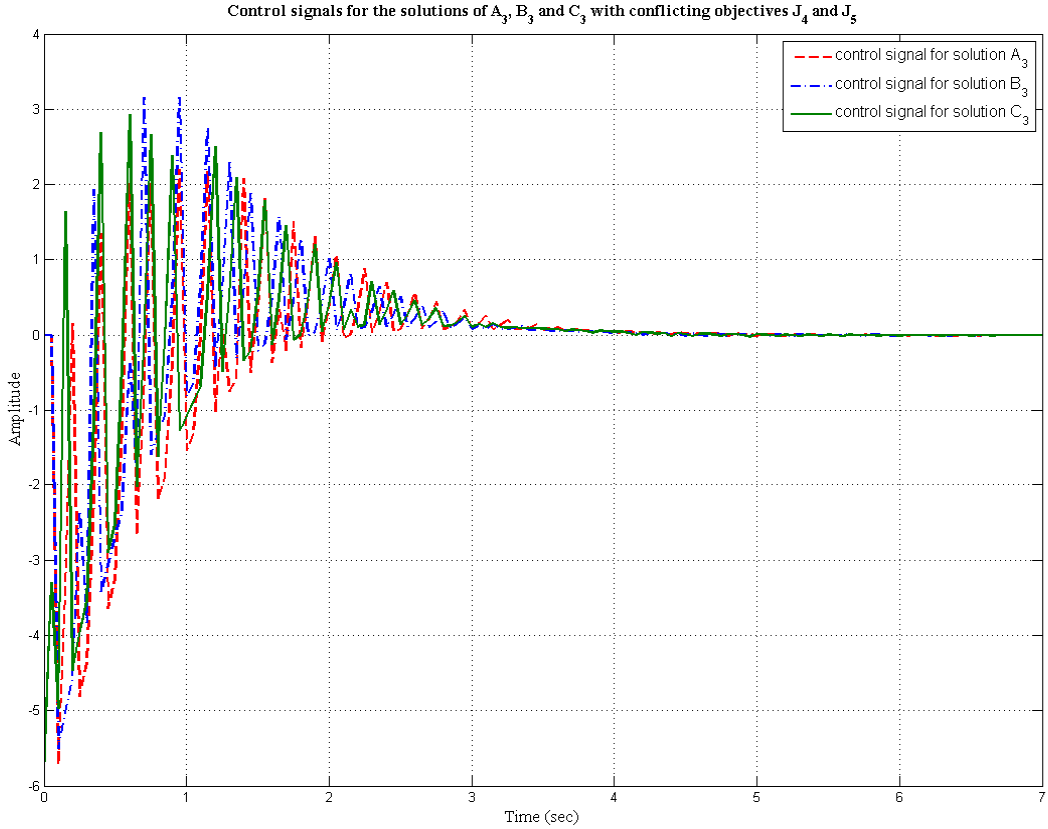

Fig. 30. Plot of control signals for the solutions $A_3$, $B_3$ and $C_3$.

### *4.4. Few discussions on the proposed networked predictive controller design scheme*

In Li *et al.* [23] system's control performance is improved by optimizing the control parameters using Estimation of Distribution Algorithm (EDA) which is a single objective optimization technique. Single objective optimization only finds best solution i.e. maximum or minimum value of a single objective. It does not give a set of solutions as a trade-off between different contradictory objectives. But in the present study, a new controller design philosophy for NCS applications has been proposed using multi-objective genetic algorithm, similar to that studied in [38], [39]. Due to the discontinuous nature of the objective functions, in the presence of random packet losses, it is not possible to frame it as classical convex optimization problem which again justifies the motivation of applying population based global optimization techniques [10]. The proposed predictive controller structure as in Pan *et al.* [6] not only ensures high degree of robust stability of the NCS but also it gives optimum time domain performance which is enforced in the present work with a MOGA based approach. Other evolutionary and swarm based algorithms could also have been employed for the present controller design problem similar to that done using particle swarm optimization [40] or genetic programming [41]. But for simplicity we restricted our study to multi-objective genetic algorithm only.

The scheme proposed in Pan *et al.* [6] is mainly for safety critical applications. In case of real time safety critical application User Datagram Protocol (UDP) is used whose transmission error checking capability is very poor. But this has the additional advantage that there are no retransmissions, unlike that of Transmission Control Protocol (TCP). This is especially important in real time control applications where the data must arrive within hard time limits and the integrity of the data is less important. Also the controller is chosen as time-driven rather than event-driven [4] as in Pan *et al.* [6]. As in the time-triggered concept the total system is synchronized with a single clock, network load will be reduced drastically. Moreover time window of every system is predefined so that every subsystem can be developed independently. In this scheme also the delayed packets are considered as dropped to reduce complexity of the LMI formulation.

Inspite of the advantages, the proposed scheme has few shortcomings. It can be seen that for a system with $n$ state variables and $M$ number of consecutive packet drops, the number of solution variables (predictive controller gains) that are to be determined is $(M+1)n$. Thus the solution technique suffers from the curse of dimensionality. But for high dimensional systems where most of the state variables cannot be measured online, a reduced order model can be obtained and then an observer can be designed based on the reduced order model. This would reduce the number of state variables and consequently reduce the number of predictive controller gains. Generally in adaptive control where the system matrices are time varying and to tackle these problems the controller gains needs to be time varying, online system identification and optimization is used. Whereas in our case the nominal system matrix is fixed and the controller gains are not time varying. But the augmented system with different packet drop situations denote different matrices which has been stabilized using a LMI based approach and then optimized to obtain the design trade-off. Since the optimization problem is offline there does not exist any hard limits on the timeline as in the case of real time online optimization [10]. Hence computational complexity and guaranteed convergence are not of major concern for this particular problem on networked predictive controller design.

It is often questioned whether such an algorithm can be implemented online or not. The proposed controller design technique is off-line. In real time applications, there is a hard time constraint and the computation must finish within a fraction of the sampling time for effective operation. Since the 'predictive' term here refers to the philosophy of applying proper control signals depending on the state of packet drop at each time step, it is a different paradigm from the traditional notions of predictive control (like MPC). Unlike MPC, this proposed methodology does not need to perform a finite horizon optimization at each time step and hence takes much less computational resources. Thus it is more suitable for real time implementation.

## 5. Conclusions

The main focus of this paper is to propose a design methodology based on multi-objective optimization by NSGA-II which gives a set of Pareto optimal solutions to optimize different contradictory time domain performance objectives for a NCS with packet drop-outs. Simulation results show the validity of the proposed approach. The multi-objective framework gives additional choices to the designer to choose the controller according to the specific requirements of the system's performance. Also chosen gains give guaranteed stability in the sense of Lyapunov since the LMI constraints are also solved as a sub-problem of the multi-objective optimization algorithm. Future works can be directed towards extending the predictive gain scheme in the event triggered systems to stabilize networked systems under random network induced delays [42], [43] along with packet drops.

In the proposed scheme we restricted the theoretical analysis for linear systems only or linearized state space models for nonlinear systems around an operating point. An optimization framework is then proposed to study design trade-offs between different time domain control objectives. The linearity of the process is assumed throughout the paper. As a sub-problem of the multi-objective optimization, the LMI formulation is introduced which guarantees Lyapunov stability inspite of arbitrary packet losses. Theoretical extension of the stability proofs for nonlinear systems with the adopted predictive state feedback law to compensate for random packet losses is thus left as the scope for future research.

**Acknowledgement:**

The authors are thankful to the anonymous reviewers for their constructive and helpful comments to enhance the quality of this paper.

**Reference:**

[1] Yodyium Tipsuwan and Mo-Yuen Chow, "Gain scheduler middleware: a methodology to enable existing controllers for networked control and teleoperation-part I: networked control", *IEEE Transactions on Industrial Electronics*, vol. 51, no. 6, pp. 1218-1227, Dec. 2004.

[2] K. C. Lee, S. Lee, and M. H. Lee, "QoS-based remote control of networked control systems via Profibus token passing protocol", *IEEE Transactions on Industrial Informatics*, vol. 1, no. 3 pp. 183-191, August 2005.

[3] John Baillieul and Panos J. Antsaklis, "Control and communication challanges in networked real-time systems", *Proceedings of the IEEE*, vol. 95, no. 1, pp. 9-28, Jan. 2007.

[4] Wei Zhang, Michael S. Branicky, and Stephen M. Phillips, "Stability of networked control systems", *IEEE Control Systems Magazine*, vol. 21, no. 1, pp. 84-99, Feb. 2001.

[5] Junxia Mu, G.P. Liu, and David Rees, "Design of robust networked predictive control systems", *Proceedings of the 2005 American Control Conference*, vol. 1, pp. 638-643, June 2005.

[6] Indranil Pan, Saptarshi Das, Soumyajit Ghosh, and Amitava Gupta, "Stabilizing gain selection of networked variable gain controller to maximize robustness using particle swarm optimization", *Proceedings of 2011 International Conference on Process Automation, Control and Computing, PACC 2011*, art no. 5978958, July 2011, Coimbatore.

[7] Abdullah Konak, David W. Coit, and Alice E. Smith, "Multi-objective optimization using genetic algorithms: A tutorial", *Reliability Engineering and System Safety*, vol. 91, no. 9, pp. 992-1007, Sept. 2006.

[8] Junlin Xiong and James Lam, "Stabilization of linear systems over networks with bounded packet loss", *Automatica*, vol. 43, no. 1, pp. 80-87, Jan. 2007.

[9] Jeremy G. VanAntwerp and Richard D. Braatz, "A tutorial on linear and bilinear matrix inequalities", *Journal of Process Control*, vol. 10, no. 4, pp. 363-385, August 2000.

[10] P.J. Fleming and R.C. Purshouse, "Evolutionary algorithms in control systems engineering: a survey", *Control Engineering Practice*, vol. 10, no. 11, pp. 1223-1241, Nov. 2002.

[11] Koan-Yuh Chang and Wen-Jer Chang, "Multi-objective control design for stochastic large-scale systems based on LMI approach and sliding mode control concept", *Journal of Marine Science and Technology*, vol. 16, no. 3, pp. 197-206, 2008.

[12] Carsten Scherer, Pascal Gahinet, and Mahmoud Chilali, "Multiobjective output-feedback control via LMI optimization", *IEEE Transactions on Automatic Control,* vol. 42, no. 7, pp. 896-911, July 1997.

[13] Alberto Herreros, Enrique Baeyens, and Jose R. Peran, "MRCD: a genetic algorithm for multiobjective robust control design", *Engineering Applications of Artificial Intelligence*, vol. 15, no. 3-4, pp. 285-301, June-August 2002.

[14] Gustavo Sanchez, Minaya Villasana, and Miguel Strefezza, "Multi-objective pole placement with evolutionary algoritms", *Evolutionary Multi-Criterion Optimization, Lecture Notes in Computer Science*, vol. 4403/2007, pp. 417-427.

[15] S. Sutha and T. Thyagarajan, "Eigenstructure assignment based multiobjective dynamic state feedback controller design for MIMO systems using NSGA-II", *The 2010 International Conference on Modelling, Identification and Control (ICMIC)*, pp. 870-875, July 2010, Okayama.

[16] Adrian Gambier, "MPC and PID control based on multi-objective optimization", *2008 American Control Conference*, pp. 4727-4732, June 2008, Seattle, WA.

[17] Adrian Gambier and Essameddin Bareddin, "Multi-objective optimal control: an overview", *IEEE Conference on Control Applications, CCA 2007*, pp. 170-175, Oct. 2007, Singapore.

[18] Adrian Gambier and Meike Jipp, "Multi-objective optimal control: an introduction", *2011 8$^{th}$ Asian Control Conference (ASCC)*, pp. 1084-1089, May 2011, Kaohsiung.

[19] Adrian Gambier, "Performance evaluation of several multi-objective optimization methods for control purposes", *2011 8$^{th}$ Asian Control Conference (ASCC)*, pp. 1067-1071, May 2011, Kaohsiung.

[20] Adrian Gambier, "Digital PID controller design based on parametric optimization", *IEEE International Conference on Control Applications, CCA 2008*, pp. 792-797, Sept. 2008, San Antonio, TX.

[21] Kalyanmoy Deb, Amrit Pratap, Sameer Agarwal, and T. Meyarivan, "A fast and elitist multiobjective genetic algorithm: NSGA-II", *IEEE Transactions on Evolutionary Computation*, vol. 6, no. 2, pp. 182-197, April 2002.

[22] M. T. Jensen, "Reducing the run-time complexity of multiobjective EAs: the NSGA-II and other algorithms", *IEEE Transactions on Evolutionary Computation*, vol. 7, no. 5, pp. 503-515, Oct. 2003.

[23] Hongbo Li, Mo-Yuen Chow, and Zengqi Sun "EDA-based speed control of a networked dc motor system with time delays and packet losses", *IEEE Transactions on Industrial Electronics*, vol. 56, no. 5, pp. 1727-1735, May 2009.

[24] Johan Lofberg, "YALMIP: a toolbox for modeling and optimization in MATLAB", *2004 IEEE International Symposium on Computer Aided Control Systems Design*, pp. 284-289, Sept. 2004, Taipei.

[25] Jian Huang, Yongji Wang, Shuang-Hua Yang, and Qi Xu, "Robust stability for remote SISO DMC controller in networked control systems", *Journal of Process Control*, vol. 19, no. 5, pp. 743-750, May 2009.

[26] Jian Guo Li, Jing Qi Yuan, and Jun Guo Lu, "Observer-based H$_\infty$ control for networked nonlinear systems with random packet losses", *ISA Transactions*, vol. 49, no. 1, pp. 39-46, Jan. 2010.

[27] Andong Liu, Li Yu, and Wen-An Zhang, "One-step receding horizon H$_\infty$ control for networked control systems with random delay and packet disordering", *ISA Transactions*, vol. 50, no. 1, pp. 44-52, Jan. 2011.

[28] Yong Zhang, Zhenxing Liu, and Bin Wang, "Robust fault detection for nonlinear networked systems with stochastic interval delay characteristics", *ISA Transactions*, vol. 50, no. 4, pp. 521-528, Oct. 2011.

[29] Mohsen Heidarinejad, Jinfeng Liu, David Munoz de la Pena, James F. Davis, and Panagiotis D. Christofides, "Handling communication disruptions in distributed model predictive control", *Journal of Process Control*, vol. 21, no. 1, pp. 173-181, Jan. 2011.

[30] Jinna Li, Qingling Zhang, Haibin Yu, and Min Cai, "Real-time guaranteed cost control of MIMO networked control systems with packet disordering", *Journal of Process Control*, vol. 21,

no. 6, pp. 967-975, July 2011.

[31] Songlin Hu and Dong Yue, "Event-triggered control design of linear networked systems with quantizations", *ISA Transactions*, vol. 51, no. 1, pp. 153-162, Jan. 2012.

[32] Hongbo Li, Mo-Yuen Chow, and Zengqi Sun "Optimal stabilizing gain selection for networked control systems with time delays and packet losses", *IEEE Transactions on Control Systems Technology*, vol. 17, no. 5, pp. 1154-1162, Sept. 2009.

[33] H. Li, M.-Y. Chow, and Z. Sun, "State feedback stabilisation of networked control systems", *IET Control Theory and Applications*, vol. 3, no. 7, pp. 929-940, July 2009.

[34] Hongbo Li, ZengQi Sun, Mo-Yuen Chow, HuaPing Liu and BaDong Chen "Stabilization of Network Control Systems with time delay and packet dropout - Part II", *2007 International Conference on Automation and Logistics*, pp. 3012-3017, Aug. 2007, Jinan.

[35] M.B.G. Cloosterman, L. Hetel, N. van de Wouw, W.P.M.H. Heemels, J. Daafouz, and H. Nijmeijer, "Controller synthesis for networked control systems", *Automatica*, vol. 46, no. 10, pp. 1584-1594, Oct. 2010.

[36] Indranil Pan, Saptarshi Das, and Amitava Gupta, "Handling packet dropouts and random delays for unstable delayed processes in NCS by optimal tuning of $PI^\lambda D^\mu$ controllers with evolutionary algorithms", *ISA Transactions*, vol. 50, no. 4, pp. 557-572, Oct. 2011.

[37] Hu Shousong and Zhu Qixin, "Stochastic optimal control and analysis of stability of networked control systems with long delay", *Automatica*, vol. 39, no. 11, pp. 1877-1884, Nov. 2003.

[38] Mohammad Nasser Saadatzi, Javad Poshtan, Mohammad Sadegh Saadatzi, and Faezeh Tafazzoli, "Novel system identification method and multi-objective-optimal multivariable disturbance observer for electric wheelchair", *ISA Transactions*, (In Press), doi: 10.1016/j.isatra.2012.06.013.

[39] David Kuhm, Dominique Knittel, and Marie-Ange Bueno, "Robust control strategies for an electric motor driven accumulator with elastic webs", *ISA Transactions*, vol. 51, no. 6, pp. 732-742, Nov. 2012.

[40] Indranil Pan, Saptarshi Das, and Amitava Gupta, "Tuning of an optimal fuzzy PID controller with stochastic algorithms for networked control systems with random time delay", *ISA Transactions*, vol. 50, no. 1, pp. 28-36, Jan. 2011.

[41] Saptarshi Das, Indranil Pan, Shantanu Das, and Amitava Gupta, "Improved model reduction and tuning of fractional-order $PI^\lambda D^\mu$ controllers for analytical rule extraction with genetic programming", *ISA Transactions*, vol. 51, no. 2, pp. 237-261, March 2012.

[42] Yongpeng Zhang, Penrose Cofie, Augustine N. Ajuzie, Jian Zhang, and Cajetan M. Akujuobi, "Real-time random delay compensation with prediction-based digital redesign", *ISA Transactions*, vol. 50, no. 2, pp. 207-212, April 2011.

[43] ZhiLe Xia, JunMin Li, and JiangRong Li, "Delay-dependent fuzzy static output feedback control for discrete-time fuzzy stochastic systems with distributed time-varying delays", *ISA Transactions*, vol. 51, no. 6, pp. 702-712, Nov. 2012.

no. 6, pp. 967-975, July 2011.

[31] Songlin Hu and Dong Yue, "Event-triggered control design of linear networked systems with quantizations", *ISA Transactions*, vol. 51, no. 1, pp. 153-162, Jan. 2012.

[32] Hongbo Li, Mo-Yuen Chow, and Zengqi Sun "Optimal stabilizing gain selection for networked control systems with time delays and packet losses", *IEEE Transactions on Control Systems Technology*, vol. 17, no. 5, pp. 1154-1162, Sept. 2009.

[33] H. Li, M.-Y. Chow, and Z. Sun, "State feedback stabilisation of networked control systems", *IET Control Theory and Applications*, vol. 3, no. 7, pp. 929-940, July 2009.

[34] Hongbo Li, ZengQi Sun, Mo-Yuen Chow, HuaPing Liu and BaDong Chen "Stabilization of Network Control Systems with time delay and packet dropout - Part II", *2007 International Conference on Automation and Logistics*, pp. 3012-3017, Aug. 2007, Jinan.

[35] M.B.G. Cloosterman, L. Hetel, N. van de Wouw, W.P.M.H. Heemels, J. Daafouz, and H. Nijmeijer, "Controller synthesis for networked control systems", *Automatica*, vol. 46, no. 10, pp. 1584-1594, Oct. 2010.

[36] Indranil Pan, Saptarshi Das, and Amitava Gupta, "Handling packet dropouts and random delays for unstable delayed processes in NCS by optimal tuning of $PI^\lambda D^\mu$ controllers with evolutionary algorithms", *ISA Transactions*, vol. 50, no. 4, pp. 557-572, Oct. 2011.

[37] Hu Shousong and Zhu Qixin, "Stochastic optimal control and analysis of stability of networked control systems with long delay", *Automatica*, vol. 39, no. 11, pp. 1877-1884, Nov. 2003.

[38] Mohammad Nasser Saadatzi, Javad Poshtan, Mohammad Sadegh Saadatzi, and Faezeh Tafazzoli, "Novel system identification method and multi-objective-optimal multivariable disturbance observer for electric wheelchair", *ISA Transactions*, (In Press), doi: 10.1016/j.isatra.2012.06.013.

[39] David Kuhm, Dominique Knittel, and Marie-Ange Bueno, "Robust control strategies for an electric motor driven accumulator with elastic webs", *ISA Transactions*, vol. 51, no. 6, pp. 732-742, Nov. 2012.

[40] Indranil Pan, Saptarshi Das, and Amitava Gupta, "Tuning of an optimal fuzzy PID controller with stochastic algorithms for networked control systems with random time delay", *ISA Transactions*, vol. 50, no. 1, pp. 28-36, Jan. 2011.

[41] Saptarshi Das, Indranil Pan, Shantanu Das, and Amitava Gupta, "Improved model reduction and tuning of fractional-order $PI^\lambda D^\mu$ controllers for analytical rule extraction with genetic programming", *ISA Transactions*, vol. 51, no. 2, pp. 237-261, March 2012.

[42] Yongpeng Zhang, Penrose Cofie, Augustine N. Ajuzie, Jian Zhang, and Cajetan M. Akujuobi, "Real-time random delay compensation with prediction-based digital redesign", *ISA Transactions*, vol. 50, no. 2, pp. 207-212, April 2011.

[43] ZhiLe Xia, JunMin Li, and JiangRong Li, "Delay-dependent fuzzy static output feedback control for discrete-time fuzzy stochastic systems with distributed time-varying delays", *ISA Transactions*, vol. 51, no. 6, pp. 702-712, Nov. 2012.